\documentclass[reqno]{amsart}
\usepackage[english]{babel}
\usepackage{graphicx,subfigure}
\usepackage{amsmath,amssymb}
\usepackage{appendix}
\usepackage{color}
\definecolor{oneblue}{rgb}{0,0.0,0.75}

\usepackage{fancyhdr}

\makeatletter
\newcommand{\sech}{\mathop{\operator@font sech}}
\newcommand{\sign}{\mathop{\operator@font sign}}
\makeatother

\newtheorem{lemma}{Lemma}[section]
\newtheorem{theorem}{Theorem}[section]

\newtheorem{remark}{Remark}[section]
\numberwithin{equation}{section}

\begin{document}

\title[]{Numerical approximation of Benjamin-Ono type models for two-dimensional internal waves}

%\author[V. A. Dougalis]{Vassilios A. Dougalis}
%\address{Mathematics Department, University of Athens, 15784
%Zographou, Greece \and Institute of Applied \& Computational
%Mathematics, FO.R.T.H., 71110 Heraklion, Greece}
%\email{doug@math.uoa.gr}

\author[A. Duran]{Angel Duran}
\address{ Applied Mathematics Department,  University of
Valladolid, 47011 Valladolid, Spain}
\email{angeldm@uva.es}

%\author[L. Saridaki]{Leetha Saridaki}
%\address{Mathematics Department, University of Athens, 15784
%Zographou, Greece \and Institute of Applied \& Computational
%Mathematics, FO.R.T.H., 71110 Heraklion, Greece}
%\email{leetha.saridaki@gmail.com}


%\dedicatory{\large \it To the memory of Vassili}

\begin{abstract}
In this paper we analyze the numerical approximation with spectral methods of two asymptotic models for the two-dimensional propagation of internal waves in a two-layer system of fluids with rigid-lid condition. The semidiscretization in space of the periodic initial-value problem with the spectral Fourier-Galerkin method is introduced and some error estimates are derived. Then the semidiscrete system is integrated in time with the Implicit Midpoint Rule and several numerical experiments illustrate the performance of the corresponding full discretization.
\end{abstract}

\maketitle

%\newpage
%\tableofcontents
\section{Introduction}
\subsection{2D asymptotic models for internal waves. The ILW and BO regimes}
The present paper is concerned with the numerical approximation of two-dimensional asymptotic models for the propagation of internal waves, \cite{Lannes2,Duchene2021}, using spectral methods. They were derived in \cite{BLS2008} as approximation to the Euler equations modelling the propagation along the interface of a two-layer system of fluids with rigid-lid condition and horizontal, flat bottom. The asymptotic models under consideration here, in the form of pde systems for the interfacial deviation and the components of the velocity of the fluid, and their level of approximation to the Euler equations, are derived under specific physical conditions, described in terms of the dimensionless parameters
\begin{eqnarray*}
\epsilon_{j}=\frac{a}{d_{j}},\quad \mu_{j}=\frac{d_{j}^{2}}{\lambda^{2}},\quad j=1,2,
\end{eqnarray*}
where $d_{j}$ denotes the depth of the $j$th layer, $j=1,2$, being  $a$ and $\lambda$, respectively, the amplitude and horizontal wavelength of the interfacial wave. We assume the so-called shallow water-small ampplitude regime, which means that the upper layer is shallow ($\mu_{1}<<1$) and the interfacial deviations are small with respect to the lower layer ($\epsilon_{2}<<1$). Under this regime, two different models, depending on the size of the rest of the parameters, are introduced.
\begin{itemize}
\item[(i)] The Intermediate Long Wave (ILW) model ($\mu\sim\epsilon^{2}\sim\epsilon_{2}<<1, \mu_{2}\sim 1$). Under this regime, the following one-parameter family of pdes are derived in \cite{BLS2008}
\begin{eqnarray}
\left[Id+\frac{\alpha\sqrt{\mu}}{\gamma}|D|{\rm coth}(\sqrt{\mu_{2}}|D|) \right]\zeta_t+\frac{1}{\gamma}\nabla\cdot\left((1- \epsilon\zeta){\bf v}\right)\nonumber&&\\
+\frac{(\alpha-1)\sqrt{\mu}}{\gamma^2}(|D|{\rm coth}(\sqrt{\mu_{2}}|D|)\nabla\cdot {\bf v}&=&0,\nonumber\\
\partial_t{\bf v}+(1-\gamma)\nabla \zeta-\frac{\epsilon}{2\gamma}\nabla|{\bf v}|^{2}&=&0.\label{BO11b}
\end{eqnarray}
In (\ref{BO11b}), $\zeta=\zeta(x,y,t)$ is the displacement of the interface at $(x,y,t)$ from some rest level, ${\bf v}=(v_{1}(x,y,t),v_{2}(x,y,t))$ represents the horizontal components of the velocity at $(x,y,t)$, and $\gamma<1$ is the density ratio between the two layers. Being $\nabla$ the gradient operator (with respect to the horizontal variables $(x,y)$) and $\Delta$ the Laplace operator, by $|D|$ we denote the operator $|D|=(-\Delta)^{1/2}$, that has the Fourier symbol
\begin{eqnarray*}
\widehat{\left(|D|f\right)}({\bf k})=|{\bf k}|\widehat{f}({\bf k}),\quad {\bf k}\in \mathbb{R}^{2},
\end{eqnarray*}
where $\widehat{f}({\bf k})$ is the Fourier transform of $f\in L^{2}(\mathbb{R}^{2})$ at ${\bf k}=(k_{x},k_{y})$ and $|{\bf k}|=\sqrt{k_{x}^{2}+k_{y}^{2}}$.
\item[(ii)] The Benjamin-Ono (BO) model ($\mu\sim \epsilon^{2}<<1, \mu_{2}=\infty$). This corresponds to the limiting case $\mu_{2}\rightarrow \infty$ in (\ref{BO11b}) (lower layer of infinite depth). The pde system becomes
\begin{eqnarray}
\left[Id+\frac{\alpha\sqrt{\mu}}{\gamma}|D|\right]\zeta_t+\frac{1}{\gamma}\nabla\cdot\left((1- \epsilon\zeta){\bf v}\right)
+\frac{(\alpha-1)\sqrt{\mu}}{\gamma^2}|D|\nabla\cdot {\bf v}&=&0,\nonumber\\
\partial_t{\bf v}+(1-\gamma)\nabla \zeta-\frac{\epsilon}{2\gamma}\nabla|{\bf v}|^{2}&=&0.\label{BO12b}
\end{eqnarray}
In (\ref{BO11b}), (\ref{BO12b}), $\alpha$ is a modelling parameter introduced to ensure well-posedness properties of the corresponding initial-value problem (ivp) with the so-called BBM trick, \cite{BBM72,S}, and by comparison with related models, \cite{CGK}. Thus, for linear well-posedness, $\alpha\geq 1$ is required, \cite{BLS2008}. In addition, nonlinear well-posedness of the ivp for (\ref{BO11b}) and  (\ref{BO12b}) is studied in \cite{X}, establishing local and long time existence/uniqueness of solutions when $\alpha> 1$.
\end{itemize}
\subsection{The periodic vip}
The main purpose of the present paper is the numerical approximation with spectral methods of the periodic ivp of (\ref{BO11b}) and (\ref{BO12b}) written in a unified formulation as
\begin{eqnarray}
\left[Id+g(D)\right]\zeta_t+\frac{1}{\gamma}\left[Id+\frac{\alpha-1}{\alpha}g(D)\right]\nabla\cdot{\bf v}-\frac{\epsilon}{\gamma}\nabla\cdot(\zeta{\bf v})&=&0,\nonumber\\
\partial_t{\bf v}+(1-\gamma)\nabla \zeta-\frac{\epsilon}{2\gamma}\nabla|{\bf v}|^{2}&=&0.\label{BO13a}
\end{eqnarray}
in $\Omega=(0,2\pi)^{2}$ for simplicity, $t>0$, and initial conditions
\begin{eqnarray}
\zeta(x,y,0)=\zeta_{0}(x,y),\quad {\bf v}(x,y,0)={\bf v}^{0}(x,y)=(v_{1}^{0}(x,y),v_{2}^{0}(x,y)),\label{BO13b}
\end{eqnarray}
for $(x,y)\in\Omega$ and $2\pi$ periodic functions $\zeta_{0},v_{j}^{0}, j=1,2$. The operator $g(D)$ is defined as
\begin{equation}  \label{BO13c}
g(D)=\left\{\begin{matrix}
\frac{\alpha\sqrt{\mu}}{\gamma}|D|&\mbox{in the {BO} case},\\
\frac{\alpha\sqrt{\mu}}{\gamma}|D|\coth(\sqrt{\mu_{2}}|D|)&\mbox{\, in the {ILW} case},  
\end{matrix}\right. \ 
\end{equation}
with Fourier symbol
\begin{equation}  \label{BO13d}
g({\bf k})=\left\{\begin{matrix}
|{\bf k}|&\mbox{in the {BO} case},\\
|{\bf k}|\coth(\sqrt{\mu_{2}}|{\bf k}|)&\mbox{\, in the {ILW} case},  
\end{matrix}\right. \ 
\end{equation}
where ${\bf k}=(k_{x},k_{y})\in\mathbb{Z}^{2}$ and in the BO case, $g(0)=\alpha/\gamma$.

The arguments used in \cite{X} to prove well-posedness of the ivp for  (\ref{BO11b}) and (\ref{BO12b}) can be applied to obtain a similar result for the periodic ivp (\ref{BO13a}), (\ref{BO13b}). Technical tools considered in \cite{X}, like commutator estimates and their application to the operator $g(D)$, are also valid in the periodic setting (via Fourier representation with Fourier series). The proof is based on the following points:
\begin{itemize}
\item A suitable change of variables (cf. section \ref{sec2}) 
\begin{eqnarray}
\eta=\zeta-\frac{1}{2\gamma(1-\gamma)}|{\bf v}|^{2},\quad {\bf u}={\bf v},\label{BO12c}
\end{eqnarray}
which satisfies, from (\ref{BO13a}), (\ref{BO13b}), the system
\begin{eqnarray}
\partial_{t}\eta-\frac{1}{\gamma}({\bf u}\cdot\nabla)\eta-\frac{1}{\gamma}(Id+g(D))^{-1}\nabla\cdot\left({\bf u}(\eta+\frac{1}{2\gamma(1-\gamma)}|{\bf u}^{2}|)\right)&&\nonumber\\
+\frac{\alpha-1}{\alpha}\nabla\cdot{\bf u}+\frac{1}{\alpha\gamma}(Id+g(D))^{-1}\nabla\cdot{\bf u}&=&0,\nonumber\\
\partial_t{\bf u}+(1-\gamma)\nabla \eta&=&0.\label{BO12d}\\
\eta(0)=\eta_{0}:=\zeta_{0}-\frac{1}{2\gamma(1-\gamma)}|{\bf v}^{0}|^{2},\quad {\bf u}(0)={\bf v}^{0}.&&\nonumber
\end{eqnarray}
\item A priori estimates.
\item Existence and uniquenes of solution of the linearized system with source terms
\begin{eqnarray*}
\partial_{t}\zeta+\frac{\alpha-1}{\alpha\gamma}\nabla\cdot{\bf v}&=&f_{1},\\
\partial_{t}{\bf v}+(1-\gamma)\nabla\zeta&=&f_{2}.
\end{eqnarray*}
\item Construction of approximate solutions from the iteration
\begin{eqnarray*}
\partial_{t}\zeta^{n+1}+\frac{\alpha-1}{\alpha\gamma}\nabla\cdot{\bf v}^{n+1}&=&\frac{1}{\gamma}(Id+g(D))^{-1}\nabla\cdot(\zeta^{n}{\bf v}^{n})\\
&&-\frac{1}{\gamma\alpha}(Id+g(D))^{-1}\nabla\cdot{\bf v}^{n},\\
\partial_{t}{\bf v}+(1-\gamma)\nabla\zeta&=&\frac{1}{\gamma}({v_{1}}^{n}\nabla v_{1}^{n+1}+v_{2}^{n}\nabla v_{2}^{n+1}),
\end{eqnarray*}
and convergence to a solution of (\ref{BO11b}). The BO case  (\ref{BO12b}) is considered as the limit $\mu_{2}\rightarrow\infty$.
\end{itemize}
In addition, consistency results of solutions of both systems are proved in  \cite{X}.

Based on the same points and adapting the tools, mentioned above, to the periodic case, we may derive a similar result of existence for the periodic ivp. The statement below makes use of the following notation. By $(\cdot,\cdot)$ and $||\cdot||$ we denote, respectively, the inner product and associated norm in $L^{2}=L^{2}(\Omega)$. For the space $L^{2}\times L^{2}$, the inner product will be denoted, for ${\bf u}, {\bf v}\in L^{2}\times L^{2}$, by
\begin{eqnarray*}
\langle {\bf u}, {\bf v}\rangle=(u_{1},v_{1})+(u_{2},v_{2}),\quad {\bf u}=(u_{1},u_{2}), {\bf v}=(v_{1},v_{2}).
\end{eqnarray*}
In addition,  $|\cdot|_{\infty}$ will stand for the norm in the Banach space $L^{\infty}=L^{\infty}(\Omega)$.

For $s\geq 0$, $H^{s}=H^{s}(\Omega)$ is the $L^{2}$-based Sobolev space on $\Omega$ of order $s$ (when $s=0$, $H^{0}:=L^{2}$) and $||\cdot||_{s}$ stands for the usual inner product. For the sake of simplicity, the norms in the product spaces $H^{s}\times H^{s}$ and $L^{\infty}\times L^{\infty}$ will be denoted with the same notation, as
\begin{eqnarray*}
||{\bf u}||_{s}=(||u_{1}||_{s}^{2}+||u_{2}||_{s}^{2})^{1/2},\; |{\bf u}|_{\infty}=\max\{|u_{1}|_{\infty},|u_{2}|_{\infty}\},\; {\bf u}=(u_{1},u_{2}).
\end{eqnarray*}

For $T>0$ let $k\geq 0$ be an integer. The space of $k$th-order continuously differentiable functions $u:[0,T]\rightarrow H^{s}$ is denoted by  $C^{k}(0,T,H^{s}), s\geq 0$ with norm
$$||u||_{T}=\max_{0\leq t\leq T}||u(s)||_{s}.$$
Finally, $f\lesssim g$ means $f\leq Cg$ for some constant $C$.
\begin{theorem}
\label{th1} Let $s>2$, $(\zeta_{0},{\bf v}^{0})\in (H^{s})^{3}$. Then there is some $T>0$, depending on $||\zeta_{0}||_{s}, ||{\bf v}^{0}||_{s}$, $\alpha$, and $\gamma$, such that (\ref{BO13a}), (\ref{BO13b}) has a unique solution $(\zeta,{\bf v})\in C^{1}(0,T,H^{s-1})\times C(0,T,H^{s})^{2}$.
\end{theorem}
\subsection{Highligths and structure of the paper}
The first purpose of the paper is the study of the numerical approximation of (\ref{BO13a}), (\ref{BO13b}) using spectral methods to discretize the spatial variables. This approach has been successfully used in \cite{DDS2} for the 1D version of (\ref{BO13a}) and in \cite{D2026_1,D2026_2} for related 2D models for internal wave propagation. These references are in the origin of the present work.

Two main reasons motivate the use of spectral methods to approximate (\ref{BO13a}). One is the presence of nonlocal terms, which can be treated directly in the spectral discretization with discrete Fourier techniques, while other approaches require the implementation of hybrid procedures, \cite{DDM}. A second reason is concerned with the accuracy of the approximation, the strongest point of the spectral discretizations, which is typically associated with the smoothness of the problem, with the more regular the solution the more accurate the approximation. This is in contrast with other approaches, \cite{DMS2007,DMS3}. In this sense, the derivation of error estimates in the present paper follows a different strategy from that used in \cite{DDS2} for the 1D case. This is mainly due to the differences between one and two dimensions as far as projection error estimates and inverse inequalities are concerned, as well as the treatment of some nonlinear terms, valid in 1D but not in 2D. The approach considered in this paper for the two-dimensional problem is based on the change of variables introduced in \cite{X} to study well-posedness, but applied to the numerical solution and the error functions. The resulting system allows to overcome the mentioned drawbacks and extend the results in 1D to 2D, although  for higher-order Sobolev spaces; that is, more regularity for the exact solution is required. On the contrary, while in 1D the error estimates are not optimal for any of the variables, the ones derived here for 2D are optimal in the $H^{2}$ norm for the ${\bf v}$ variables, and suboptimal for the $\zeta$ component. The existence of numerical solution and error analysis are the contents of section \ref{sec2}.

The second purpose of this paper is to complete the approximation of (\ref{BO13a}), (\ref{BO13b}) with the introduction of a fully discrete scheme from the semisdiscretization given by the spectral approach. The choice of a standard time integrator as the Implicit Midpoint Rule (IMR) is due to several reasons. The main are concerned with the method itself, its beneficial stability, geometric, and implementation properties for this kind of problems, \cite{HLW2004}. They are particularly useful in the model considered here, where the semidiscrete system is stiff and the investigation of the model may require long-time simulations. In addition, although the method is implicit, it shows some implementation advantages in the iterative resolution of the internal stages when applied to the semidiscrete system. The formulation and performance of the fully discrete system are analyzed in section \ref{sec3} and the numerical experiments introduced there can be compared with the results obtained by similar approximations to other 2D models for internal waves in \cite{D2026_1,D2026_2}.

\section{Spectral discretization}
\label{sec2}
In this section we introduce and analyze the semidiscrete Fourier Galerkin spectral approximation of the periodic ivp (\ref{BO13a}), (\ref{BO13b}), written here in the form
\begin{eqnarray}
\partial_{t}\zeta +\frac{1}{\gamma}\mathcal{J}_{\alpha}\nabla\cdot{\bf v}-\frac{1}{\gamma}\mathcal{T}_{g}\nabla\cdot(\zeta{\bf v})&=&0,\nonumber\\
\partial_t{\bf v}+(1-\gamma)\nabla \zeta-\frac{\epsilon}{2\gamma}\nabla|{\bf v}|^{2}&=&0.\label{BO21}\\
\zeta(0)=\zeta_{0},\quad {\bf v}(0)={\bf v}^{0},&&\nonumber
\end{eqnarray}
for $\alpha>1$. In (\ref{BO21}), two nonlocal operators are introduced:
\begin{eqnarray*}
\mathcal{T}{g}=(Id+g(D))^{-1},\quad \mathcal{J}_{\alpha}=\mathcal{T}_{g}(Id+\frac{\alpha-1}{\alpha}g(D)),\label{BO21a}
\end{eqnarray*}
where $g(D)$ is given in (\ref{BO13c}). The corresponding Fourier symbols, cf. (\ref{BO13d}), are, for ${\bf k}=(k_{x},k_{y})^{T}\in\mathbb{Z}^{2}$
\begin{eqnarray}
&&\widehat{\mathcal{T}_{g}f}({\bf k})=\tau({\bf k})\widehat{f}({\bf k}),\quad \tau({\bf k})=\frac{1}{1+g({\bf k})},\nonumber\\
&&\widehat{\mathcal{J}_{\alpha}f}({\bf k})=j_{\alpha}({\bf k})\widehat{f}({\bf k}),\quad j_{\alpha}({\bf k})=(1+\left(\frac{\alpha-1}{\alpha}\right)g({\bf k}))\tau({\bf k}).\label{BO22}
\end{eqnarray}
The following properties of $g({\bf k})$, in the ILW case, derived in \cite{DDS2}, will be used in the sequel:
\begin{eqnarray}
g({\bf k})&\approx&\frac{\alpha}{\gamma}|{\bf k}|,\quad |{\bf k}|\geq 1,\nonumber\\
g({\bf k})&\approx&\frac{\alpha}{\gamma}\left(1+\frac{|{\bf k}|^{2}}{3}+O(|{\bf k}|^{4})\right),\quad |{\bf k}|\; {\rm small}.\label{BO23}
\end{eqnarray}
Then, from (\ref{BO23}), it transpires that the symbols $\tau({\bf k})$ and $j_{\alpha}({\bf k})$ in (\ref{BO22}) are bounded for all ${\bf k}$ and, in addition, if $|{\bf k}|\rightarrow\infty$, then $\tau({\bf k})=O\left(\frac{1}{|{\bf k}|}\right)$. These properties will be useful in the error analysis made in section \ref{sec22}.
\subsection{Spectral semidiscretization in space}
\label{sec21}
We consider, for an integer $N\geq 1$ and for $T>0$, the finite dimensional space
\begin{eqnarray*}
S_{N}={\rm span}\{e^{i{\bf x}\cdot{\bf k}},{\bf x}=(x,y)\in \Omega, {\bf k}=(k_{x},k_{y}), k_{j}\in\mathbb{Z}, -N\leq k_{j}\leq N, j=x,y\},
\end{eqnarray*}
where ${\bf x}\cdot{\bf k}=k_{x}x+k_{y}y$,
and define the semidiscrete Fourier-Galerkin approximation of (\ref{BO13a}), (\ref{BO13b}) for $0\leq t\leq T$ as the mappings  $\zeta^{N}, v_{j}^{N}:[0,T]\rightarrow S_{N}, j=1,2,$ such that if ${\bf v}^{N}=(v_{1}^{N},v_{2}^{N})$, for $0\leq t\leq T$ 
\begin{eqnarray}
\partial_{t}\zeta^{N} +\frac{1}{\gamma}\mathcal{J}_{\alpha}\nabla\cdot{\bf v}^{N}-\frac{1}{\gamma}\mathcal{T}_{g}\nabla\cdot\widetilde{P}_{N}(\zeta^{N}{\bf v}^{N})&=&0,\nonumber\\
\partial_t{\bf v}^{N}+(1-\gamma)\nabla \zeta^{N}-\frac{\epsilon}{2\gamma}\nabla P_{N}\left(|{\bf v}^{N}|^{2}\right)&=&0,\label{BO24}\\
\zeta(0)=P_{N}(\zeta_{0}),\quad {\bf v}(0)=\widetilde{P}_{N}({\bf v}^{0}),&&\nonumber
\end{eqnarray}
where $P_{N}$ denotes the usual $L^{2}$-orthogonal projection onto $S_{N}$ (which commutes with any spatial partial derivative, \cite{CHQZ}, and with the operator $\mathcal{T}_{g}$) and if ${\bf u}=(u_{1},u_{2})^{T}$ then $\widetilde{P}_{N}({\bf u})=(P_{N}u_{1},P_{N}u_{2})^{T}$. In terms of the corresponding Fourier coefficients $\widehat{\zeta}^{N}=\widehat{\zeta}^{N}({\bf k},\cdot), \widehat{v_{j}^{N}}=\widehat{v_{j}^{N}}({\bf k},\cdot), j=1,2, {\bf k}=(k_{x},k_{y})^{T}\in\mathbb{Z}^{2}$, the system (\ref{BO24}) reads
\begin{eqnarray}
\partial_{t}\widehat{\zeta}^{N}+\frac{1}{\gamma}j_{\alpha}({\bf k})(ik_{x}\widehat{v_{1}^{N}}+ik_{y}\widehat{v_{2}^{N}})-\tau({\bf k})({\bf k})(ik_{x}\widehat{\zeta^{N}v_{1}^{N}}+ik_{y}\widehat{\zeta^{N}v_{2}^{N}})&=&0,\nonumber\\
\partial_{t}\widehat{{\bf v}^{N}}+(1-\gamma)i{\bf k}\widehat{\zeta^{N}}-\frac{1}{2\gamma}i{\bf k}\widehat{|{\bf v}^{N}|^{2}}&=&0,\label{BO25}\\
\widehat{\zeta^{N}}(0)=\widehat{\zeta_{0}}({\bf k}),\quad \widehat{{\bf v}^{N}}(0)=\widehat{{\bf v}^{0}}({\bf k}),&&\nonumber
\end{eqnarray}
where if ${\bf u}=(u_{1},u_{2})^{T}$ then $\widehat{{\bf u}}({\bf k})=(\widehat{u_{1}}({\bf k}),\widehat{u_{2}}({\bf k}))^{T}$. In the following result, the existence of solution of (\ref{BO24}) is discussed.
\begin{theorem}
\label{th2}
Assume that $(\zeta_{0},{\bf v}^{0})\in (H^{2})^{3}$. Then there exist $T>0$ and a constant $C$, depending both on $||\zeta_{0}||_{2}, ||{\bf v}^{0}||_{2}, \gamma, \alpha$, and $N$, such that (\ref{BO24}) admits a unique solution $(\zeta^{N}(t), {\bf v}^{N}(t)), 0\leq t\leq T$, satisfying
\begin{eqnarray}
\max_{0\leq t\leq T}\left(||\zeta^{N}(t)||_{2}+||{\bf v}^{N}(t)||_{2}\right)\leq C.\label{BO25b}
\end{eqnarray}
\end{theorem}
\begin{proof}
Standard ode theory can be applied to (\ref{BO25}) in order to derive existence and uniqueness of solution $(\zeta^{N}(t), {\bf v}^{N}(t))$ of (\ref{BO24}), locally in time on some interval $0\leq t\leq t^{*}$. Following \cite{X}, we consider the new variables (cf. (\ref{BO12c}))
\begin{eqnarray}
\eta^{N}:=\zeta^{N}-\frac{1}{2\gamma(1-\gamma)}P_{N}(|{\bf v}^{N}|^{2}),\; {\bf u}^{N}:={\bf v}^{N},\label{BO26}
\end{eqnarray}
which, from (\ref{BO24}) and using commutativity of $P_{N}$ with partial derivatives, satisfy for $0\leq t\leq t^{*}$ (cf. (\ref{BO12d}))
\begin{eqnarray}
\partial_{t}\eta^{N}-\frac{1}{\gamma}P_{N}\left(({\bf u}^{N}\cdot\nabla)\eta^{N}\right)-\frac{1}{\gamma}\mathcal{T}_{g}\nabla\cdot\widetilde{P}_{N}\left({\bf u}^{N}(\eta^{N}+\frac{1}{2\gamma(1-\gamma)}P_{N}\left(|{\bf u}^{N}|^{2}\right))\right)&&\label{BO27}\\
+\frac{\alpha-1}{\alpha}\nabla\cdot{\bf u}^{N}+\frac{1}{\alpha\gamma}\mathcal{T}_{g}\nabla\cdot{\bf u}^{N}&=&0,\nonumber\\
\partial_t{\bf u}^{N}+(1-\gamma)\nabla \eta^{N}&=&0,\nonumber\\
\eta^{N}(0)=P_{N}(\eta_{0}),\quad {\bf u}^{N}(0)=\widetilde{P}_{N}({\bf v}^{0}).&&\nonumber
\end{eqnarray}
Let 
\begin{eqnarray}
Y_{0}:=\frac{\alpha\gamma}{\alpha-1}||\eta_{0}||_{2}^{2}+\frac{1}{1-\gamma}||{\bf u}^{0}||_{2}^{2},\quad
C_{\infty}:=(Y_{0}+\sqrt{Y_{0}+Y_{0}^{2}})^{1/2}.\label{BO27a}
\end{eqnarray}
Note that, by Sobolev embedding
\begin{eqnarray*}
|{\bf u}^{N}(0)|_{\infty}<C_{\infty}.
\end{eqnarray*}
Let $t_{N}\in (0,t^{*}]$ be the maximal time such that
\begin{eqnarray}
|{\bf u}^{N}(t)|_{\infty}\leq 2C_{\infty},\quad 0<t\leq t_{N}.\label{BO28}
\end{eqnarray}
This condition will be used in the following derivation of estimates for $(\eta^{N},{\bf u}^{N})$. If $x_{i},x_{j}=x,y$, we multiply the first equation of (\ref{BO27}) by
\begin{eqnarray*}
E_{ij}^{N}=\frac{\alpha\gamma}{\alpha-1}\left(\eta^{N}-\partial_{x_{i}}^{2}\eta^{N}+\partial_{x_{i}}^{2}\partial_{x_{j}}^{2}\eta^{N}\right),
\end{eqnarray*}
the second equations by
\begin{eqnarray*}
U_{ij}^{N}=\frac{1}{1-\gamma}\left({\bf u}^{N}-\partial_{x_{i}}^{2}{\bf u}^{N}+\partial_{x_{i}}^{2}\partial_{x_{j}}^{2}{\bf u}^{N}\right),
\end{eqnarray*}
(componentwise), integrate over $\Omega$, and add the resulting equations, \cite{X}. Integrating by parts, using periodicity, commutativity of $P_{N}$ with $\mathcal{T}_{g}$ and the partial derivatives $\partial_{x_{i}}$, and the property
\begin{eqnarray*}
(P_{N}f,\varphi)=(f,\varphi),\quad f\in L^{2}, \varphi\in S_{N},
\end{eqnarray*}
we have, for $0<t\leq t_{N}$
\begin{eqnarray}
\frac{d}{dt}\left(\frac{\alpha\gamma}{\alpha-1}||\eta^{N}(t)||_{2}^{2}+\frac{1}{1-\gamma}||{\bf u}^{N}(t)||_{2}^{2}\right)=A_{1}+A_{2}+A_{3}\label{BO29a},
\end{eqnarray}
where
\begin{eqnarray}
A_{1}&=&
{-\sum_{i,j}\int_{\Omega}\frac{1}{\alpha-1}(\mathcal{T}_{g}\nabla\cdot{\bf u}^{N})E_{ij}^{N}d\underline{x}},\nonumber\\
A_{2}&=&{\sum_{i,j}\int_{\Omega}\frac{\alpha}{\alpha-1}\left(({\bf u}^{N}\cdot\nabla)\eta^{N}\right)E_{ij}^{N}d\underline{x}},\nonumber\\
A_{3}&=&{\sum_{i,j}\int_{\Omega}\frac{\alpha}{\alpha-1}\mathcal{T}_{g}\nabla\cdot\left({\bf u}^{N}(\eta^{N}+\frac{1}{2\gamma(1-\gamma)}|{\bf u}^{N}|^{2})\right)E_{ij}^{N}d\underline{x}}.\label{BO29b}
\end{eqnarray}
We now estimate the right-hand side of (\ref{BO29a}) from the analysis of each term in  (\ref{BO29b}). Using commutativity between $\mathcal{T}_{g}\nabla$ and partial derivatives $\partial_{x_{i}}$ and since, from (\ref{BO23}), $\mathcal{T}_{g}\nabla$ is bounded, then we have, for $0<t\leq t_{N}$
\begin{eqnarray}
|A_{1}|\lesssim ||{\bf u}^{N}||_{2}||\eta^{N}||_{2}.\label{BO210}
\end{eqnarray}
Similar arguments to those made in \cite{X} can apply to write $A_{2}$, after integration by parts, as
\begin{eqnarray*}
A_{2}&=&\frac{\alpha}{\alpha-1}\sum_{i,j}\left(\int_{\Omega}\left(({\bf u}^{N}\cdot\nabla)\eta^{N}\right)\eta^{N}+\partial_{x_{i}}\left(({\bf u}^{N}\cdot\nabla)\eta^{N}\right)\partial_{x_{i}}\eta^{N}\right.\nonumber\\
&&\left.+\partial_{x_{i}}\partial_{x_{j}}\left(({\bf u}^{N}\cdot\nabla)\eta^{N}\right)\partial_{x_{i}}\partial_{x_{j}}\eta^{N}\right)d\underline{x}\nonumber\\
&=&\frac{\alpha}{\alpha-1}\sum_{i,j}\left(\int_{\Omega}\left(({\bf u}^{N}\cdot\nabla)\eta^{N}\right)\eta^{N}+[\partial_{x_{i}},{\bf u}^{N}]\cdot\nabla\eta^{N}\partial_{x_{i}}\eta^{N}\right.\\
&&\left.+[\partial_{x_{i}}\partial_{x_{j}},{\bf u}^{N}]\cdot\nabla\eta^{N}\partial_{x_{i}}\partial_{x_{j}}\eta^{N}-\frac{1}{2}\left(\nabla\cdot{\bf u}^{N}\right)\left(\partial_{x_{i}}\eta^{N}\right)^{2}\right.\\
&&\left.-\frac{1}{2}\left(\nabla\cdot{\bf u}^{N}\right)\left(\partial_{x_{i}}\partial_{x_{j}}\eta^{N}\right)^{2}\right)d\underline{x},
\end{eqnarray*}
where the commutator 
\begin{eqnarray}
[L,f]g=L(fg)-fLg,\label{BO210a}
\end{eqnarray}
is used by components, that is
\begin{eqnarray*}
[L,{\bf u}^{N}]\cdot\nabla\eta^{N}=[L,u_{1}^{N}]\partial_{x}\eta^{N}+[L,u_{2}^{N}]\partial_{y}\eta^{N},
\end{eqnarray*}
where $L=\partial_{x_{i}}$ or $\partial_{x_{i}}\partial_{x_{j}}$. The definition of the commutator operator (\ref{BO210a}) and Cauchy-Schwarz inequality imply that
\begin{eqnarray}
|A_{2}|\lesssim ||{\bf u}^{N}||_{2}\left(||\eta^{N}||_{2}+||\eta^{N}||_{2}^{2}\right).\label{BO211}
\end{eqnarray}
We finally decompose $A_{3}=A_{31}+A_{32}$, where
\begin{eqnarray*}
A_{31}&=&{\sum_{i,j}\int_{\Omega}\frac{\alpha}{\alpha-1}\mathcal{T}_{g}\nabla\cdot\left({\bf u}^{N}(\eta^{N})\right)E_{ij}^{N}d\underline{x}}
\\
A_{32}&=&{\sum_{i,j}\int_{\Omega}\frac{\alpha}{\alpha-1}\frac{1}{2\gamma(1-\gamma)}\mathcal{T}_{g}\nabla\cdot\left({\bf u}^{N}(|{\bf u}^{N}|^{2})\right)E_{ij}^{N}d\underline{x}}.
\end{eqnarray*}
From boundedness of $\mathcal{T}_{g}\nabla$ and commutativity with the $\partial_{x_{i}}$, we have, after integration by parts
\begin{eqnarray}
|A_{31}|\lesssim ||{\bf u}^{N}||_{2}||\eta^{N}||_{2}^{2}.\label{BO212}
\end{eqnarray}
From the same properties and using (\ref{BO28}), we have
\begin{eqnarray}
|A_{32}|\lesssim ||{\bf u}^{N}||_{2}^{2}||\eta^{N}||_{2}.\label{BO213}
\end{eqnarray}
Therefore, from (\ref{BO210}), (\ref{BO211})-(\ref{BO213}), applied to (\ref{BO29a}), we obtain the estimate
\begin{eqnarray}
\frac{d}{dt}\left(\frac{\alpha\gamma}{\alpha-1}||\eta^{N}||_{2}^{2}+\frac{1}{1-\gamma}||{\bf u}^{N}||_{2}^{2}\right)&\leq& C\left(||{\bf u}^{N}||_{2}^{2}+||\eta^{N}||_{2}^{2}+||{\bf u}^{N}||_{2}||\eta^{N}||_{2}^{2}\right.\nonumber\\
&&\left.+||{\bf u}^{N}||_{2}^{2}||\eta^{N}||_{2}\right),\label{BO214}
\end{eqnarray}
for some constant $C$ and $0<t\leq t_{N}$. We now consider, cf. \cite{X}
\begin{eqnarray*}
Y(t):=\frac{\alpha\gamma}{\alpha-1}||\eta^{N}(t)||_{2}^{2}+\frac{1}{1-\gamma}||{\bf u}^{N}(t)||_{2}^{2},\quad 0<t\leq t_{N},
\end{eqnarray*}
where $Y(0)=Y_{0}$, defined in (\ref{BO27a}). Thus, (\ref{BO214}) implies
\begin{eqnarray*}
\frac{d}{dt}Y(t)\lesssim (Y(t)+Y(t)^{3/2})\leq CY(t)(1+Y(t)),
\end{eqnarray*}
for some constant $C$ and $0<t\leq t_{N}$. Then, defining
\begin{eqnarray*}
T=\frac{1}{2C}\ln\left(1+\frac{1}{Y_{0}}\right),
\end{eqnarray*}
we have, for $0<t\leq T$, cf. \cite{X}
\begin{eqnarray*}
Y(t)&\leq & \frac{1}{\left(1+\frac{1}{Y_{0}}\right)e^{-tC}-1}\leq \frac{1}{\left(1+\frac{1}{Y_{0}}\right)^{1/2}-1}\\
&=&\left(1+\left(1+\frac{1}{Y_{0}}\right)^{1/2}\right)Y_{0}=(Y_{0}+\sqrt{Y_{0}+Y_{0}^{2}})^{1/2}=C_{\infty}.
\end{eqnarray*}
Hence $T\leq t_{N}$ and using (\ref{BO26}) and the embedding $H^{2}\hookrightarrow L^{\infty}$, then (\ref{BO25b}) holds.
\end{proof}
\subsection{Error analysis}
\label{sec22}
In order to derive some error estimates for (\ref{BO24}), two main tools will be use. One concerns standard projection error estimates in 2D, inverse inequalities in $S_{N}$ and Grisvard's lemma. More specifically:
\begin{itemize}
\item[(i)] Projection errors, \cite{CHQZ}: { For integers $0\leq j\leq \mu$ and $u\in H^{\mu}(\Omega)$
\begin{eqnarray}
||u-P_{N}u||_{{j}}&\lesssim &N^{j-\mu}||u||_{{\mu}},\label{BO215}\\
|u-P_{N}u|_{\infty}&\lesssim & N^{1-\mu}||u||_{{\mu}},\quad \mu> 1.\nonumber
\end{eqnarray}
}
\item[(ii)] Inverse inequalities, \cite{CHQZ}: {For integers $0\leq l\leq m$, $\phi\in S_{N}$
\begin{eqnarray}
||\phi||_{{m}}&\lesssim &N^{m-l}||\phi||_{{l}},\label{BO216}\\
||\phi||_{m,\infty}&\lesssim &N^{m+1-l}||\phi||_{{l}}.\nonumber
\end{eqnarray}
}
\item[(iii)] A proof of Grisvard's lemma can be found in \cite{Grisvard}:
\begin{lemma}
Let $s_{j}\in\mathbb{R}, j=1,2,3$ be such that $s_{1},s_{2}\geq s_{3}, s_{1}+s_{2}\geq , s_{1}+s_{2}-s_{3}>1$. Then, $(f_{1},f_{2})\mapsto f_{1}f_{2}$ is bilinear continuous from $H^{s_{1}}\times H^{s_{2}}$ to $H^{s_{3}}$.
\end{lemma}
\end{itemize}
The second tool makes use of the change of variables (\ref{BO26}) as well as the corresponding (\ref{BO12c}) for the continuous problem (\ref{BO13a}). As usual, we define
\begin{eqnarray*}
\theta=\zeta^{N}-P_{N}\zeta,&& \rho=P_{N}\zeta-\zeta,\\
{\bf \xi}=(\xi_{1},\xi_{2})={\bf u}^{N}-\widetilde{P_{N}}{\bf u},&& {\bf \sigma}=(\sigma_{1},\sigma_{2})=\widetilde{P_{N}}{\bf u}-{\bf u},
\end{eqnarray*}
so that $\zeta^{N}-\zeta=\theta+\rho$, ${\bf u}^{N}-{\bf u}={\bf \xi}+{\bf \sigma}$. By comparing (\ref{BO12d}) and (\ref{BO27}), the following problem for $(\theta,\xi)$ holds
\begin{eqnarray}
\partial_{t}\theta+\frac{\alpha-1}{\alpha}\nabla\cdot\xi+\frac{1}{\alpha\gamma}\mathcal{T}_{g}\nabla\cdot\xi&=&\frac{1}{\gamma}P_{N}\widetilde{A},\nonumber\\
\partial_{t}\xi+(1-\gamma)\nabla \theta&=&0,\label{BO220}\\
\theta(0)=0,\quad \xi(0)=0.&&\nonumber
\end{eqnarray}
where
\begin{eqnarray*}
\widetilde{A}&=&({\bf u}\cdot\nabla)\eta-P_{N}\left(({\bf u}^{N}\cdot\nabla)\eta^{N}\right)+\mathcal{T}_{g}\nabla\cdot\left({\bf u}(\eta+\frac{1}{2\gamma(1-\gamma)}|{\bf u}|^{2})\right.\nonumber\\
&&\left.-\widetilde{P}_{N}\left({\bf u}^{N}(\eta^{N}+\frac{1}{2\gamma(1-\gamma)}|{\bf u}^{N}|^{2})\right)\right).\label{BO221}
\end{eqnarray*}
The key result is the derivation of error estimates for the transformed systems. Using Theorems \ref{th1} and \ref{th2}, let $T_{0}>0$ and $T_{1}>0$ such that (\ref{BO12d}) admits a unique solution for $0<t\leq T_{0}$ and (\ref{BO27}) admits a unique solution for $0<t\leq T_{1}$.
\begin{theorem}
\label{th3} Let $\alpha>1$ and assume that the solution $(\eta,{\bf u})$ of (\ref{BO12d}) satisfies $\eta(t)\in H^{\mu_{1}}, \mu_{1}>3, {\bf u}(t)\in (H^{\mu_{2}})^{2}, \mu_{2}>2, 0\leq t\leq T_{0}$. Then, for $N$ large enough, there are some $T^{*}=T^{*}(N)>0$ and a constant $C$, independent of $N$, such that
\begin{eqnarray}
\max_{0\leq t\leq T^{*}}\left(||\eta^{N}(t)-\eta(t)||_{2}+||{\bf u}^{N}(t)-{\bf u}(t)||_{2}\right)\leq C\left(N^{3-\mu_{1}}+N^{2-\mu_{2}}\right).\label{BO222}
\end{eqnarray}
\end{theorem}
\begin{proof}
Let $T_{2}=\min\{T_{0},T_{1}\}$, so that $\theta$ and $\xi$ are defined and satisfy (\ref{BO220}) for $0\leq t\leq T_{2}$. Since $\xi(0)=0$, there is a maximal value of time $t_{N}>0$ such that
\begin{eqnarray}
|\xi(t)|_{\infty}\leq 1,\quad 0\leq t\leq t_{N}.\label{maximalt}
\end{eqnarray}

In order to derive (\ref{BO222}), we will follow a similar strategy to that used in the proof of Theorem \ref{th2}. For $x_{i},x_{j}=x,y$, we multiply the first equation of (\ref{BO220}) by
\begin{eqnarray*}
\Theta_{ij}^{N}=\frac{\alpha\gamma}{\alpha-1}\left(\theta-\partial_{x_{i}}^{2}\theta+\partial_{x_{i}}^{2}\partial_{x_{j}}^{2}\theta\right),
\end{eqnarray*}
the second equations by
\begin{eqnarray*}
\Xi_{ij}^{N}=\frac{1}{1-\gamma}\left(\xi-\partial_{x_{i}}^{2}\xi+\partial_{x_{i}}^{2}\partial_{x_{j}}^{2}\xi\right),
\end{eqnarray*}
integrate over $\Omega$, and add the resulting equations. The same arguments lead for $x_{i},x_{j}=x,y$ and $0<t\leq t_{N}$, to
\begin{eqnarray}
\frac{d}{dt}\left(\frac{\alpha\gamma}{\alpha-1}||\eta^{N}||_{2}^{2}+\frac{1}{1-\gamma}||{\bf u}^{N}||_{2}^{2}\right)=B_{1}+B_{2}\label{BO222a},
\end{eqnarray}
where
\begin{eqnarray*}
B_{1}&=&
-\sum_{i,j}\int_{\Omega}\frac{1}{\alpha-1}\left((\mathcal{T}_{g}\nabla\cdot\xi)\theta
+(\mathcal{T}_{g}\nabla\cdot\partial_{x_{i}}\xi)\partial_{x_{i}}\theta\right.\\
&&\left.
-(\mathcal{T}_{g}\nabla\cdot\partial_{x_{i}}\partial_{x_{j}}\xi)\partial_{x_{i}}\partial_{x_{j}}\right)\theta
 d\underline{x},\nonumber\\
B_{2}&=&{\sum_{i,j}\int_{\Omega}\frac{\alpha}{\alpha-1}\widetilde{A}\Theta_{ij}^{N}d\underline{x}}.
\end{eqnarray*}
Note that, since $\mathcal{T}_{g}\nabla$ is bounded, then we have, for $0<t\leq t_{N}$
\begin{eqnarray}
|B_{1}|\lesssim ||\xi||_{2}||\theta||_{2}.\label{BO222b}
\end{eqnarray}
The hardest and longest part is the estimate for $B_{2}$. We write $\widetilde{A}$ as
 $\widetilde{A}= -\widetilde{A}_{1}+ \widetilde{A}_{2}$, with
\begin{eqnarray*}
\widetilde{A}_{1}&=&\widetilde{A}_{11}+\widetilde{A}_{12}+\widetilde{A}_{13};\\
\widetilde{A}_{11}&=&{\bf u}\cdot\nabla\rho+\sigma\cdot\nabla\eta+\sigma\cdot\nabla\rho,\\
\widetilde{A}_{12}&=&\xi\cdot\nabla\eta+\xi\cdot\nabla\rho,\\
\widetilde{A}_{13}&=&{\bf u}\cdot\nabla\theta+\sigma\cdot\nabla\theta+\xi\cdot\nabla\theta,\\
\widetilde{A}_{2}&=&\mathcal{T}_{g}\nabla\cdot\left(\widetilde{A}_{21}+\widetilde{A}_{22}+\widetilde{A}_{23}+\widetilde{A}_{24}\right);\\
\widetilde{A}_{21}&=&{\bf u}\eta-({\bf u}+\xi+\sigma)(\eta+\theta+\rho)\\
&=&-({\bf u}\theta+{\bf u}\rho+\xi\eta+\xi\theta+\xi\rho+\sigma\eta+\sigma\theta+\sigma\rho),\\
\widetilde{A}_{22}&=&({\bf u}(|{\bf u}|^{2}-P_{N}(|{\bf u}|^{2}))),\\
\widetilde{A}_{23}&=&{\bf u}P_{N}(2(\xi+\sigma)\cdot {\bf u}+|\xi+\sigma|^{2})\\
&=&{\bf u}P_{N}(2\xi\cdot{\bf u}+2\sigma\cdot{\bf u}+2\xi\cdot\sigma+|\xi|^{2}+|\sigma|^{2}),\\
\widetilde{A}_{24}&=&(\xi+\sigma)P_{N}(|{\bf u}|^{2}+2{\bf u}\cdot(\xi+\sigma)+|\xi|^{2}+|\sigma|^{2}+2\sigma\cdot\xi).
\end{eqnarray*}
We now estimate each term, for $0<t\leq t_{N}$, by using arguments from \cite{X} (already considered in the proof of Theorem \ref{th2}), along with projection error estimates (\ref{BO215}) and the assumptions on the regularity of the exact solution $(\eta, {\bf u})$. For the first term of $\widetilde{A}_{11}$, we have 
\begin{eqnarray*}
({\bf u}\cdot\nabla\rho,\Theta_{ij}^{N})&=&\int_{\Omega}\left({\bf u}\cdot\nabla\rho\theta+[\partial_{x_{i}},{\bf u}]\cdot\rho\partial_{x_{i}}\theta+[\partial_{x_{i}}\partial_{x_{j}},{\bf u}]\cdot\rho\partial_{x_{i}}\partial_{x_{j}}\theta\right.\\
&&\left. +({\bf u}\cdot\nabla)\partial_{x_{i}}\rho\partial_{x_{i}}\theta+({\bf u}\cdot\nabla)\partial_{x_{i}}\partial_{x_{j}}\rho\partial_{x_{i}}\partial_{x_{j}}\theta\right)d\underline{x},
\end{eqnarray*}
and similarly for the other two terms. This implies a bound for $\widetilde{A}_{11}$ of the form
\begin{eqnarray*}
|\sum_{i,j}(\widetilde{A}_{11},\Theta_{ij}^{N})|\lesssim \left(N^{3-\mu_{1}}+N^{2-\mu_{2}}\right)||\theta||_{2}.
\end{eqnarray*}
Using the same tools, it is not hard to derive the estimates
\begin{eqnarray*}
|\sum_{i,j}(\widetilde{A}_{12},\Theta_{ij}^{N})|&\lesssim &\left(N^{3-\mu_{1}}+1\right)||\xi||_{2}||\theta||_{2},\\
|\sum_{i,j}(\widetilde{A}_{13},\Theta_{ij}^{N})|&\lesssim &\left(N^{2-\mu_{2}}+1\right)||\theta||_{2}^{2}+||\xi||_{2}||\theta||_{2}^{2}.
\end{eqnarray*}
We now study $\widetilde{A}_{2}$, for which the boundedness of $\mathcal{T}_{g}\nabla$ will be used. As for $\widetilde{A}_{21}$, observe again that
\begin{eqnarray*}
(\mathcal{T}_{g}\nabla\cdot(\theta({\bf u}+\sigma)),\Theta_{ij}^{N})&=&\int_{\Omega}\left(\mathcal{T}_{g}\nabla\cdot(\theta({\bf u}+\sigma))\theta+\mathcal{T}_{g}\nabla\cdot\partial_{x_{i}}(\theta({\bf u}+\sigma))\partial_{x_{i}}\theta\right.\\
&&\left. +\mathcal{T}_{g}\nabla\cdot\partial_{x_{i}}\partial_{x_{j}}(\theta({\bf u}+\sigma))\partial_{x_{i}}\partial_{x_{j}}\theta\right)d\underline{x},
\end{eqnarray*}
and therefore
\begin{eqnarray*}
|\sum_{i,j}(\mathcal{T}_{g}\nabla\cdot(\theta({\bf u}+\sigma)),\Theta_{ij}^{N})|\lesssim \left(N^{2-\mu_{2}}+1\right)||\theta||_{2}^{2}.
\end{eqnarray*}
A similar strategy for the other terms leads to
\begin{eqnarray*}
%|\sum_{i,j}(\widetilde{A}_{12},\Theta_{ij}^{N})|&\lesssim &\left(N^{3-\mu_{1}}+1\right)||\xi||_{2}||\theta||_{2},\\
|\sum_{i,j}(\mathcal{T}_{g}\nabla\cdot\widetilde{A}_{21},\Theta_{ij}^{N})|&\lesssim &\left(N^{2-\mu_{2}}+1\right)||\theta||_{2}^{2}+||\xi||_{2}||\theta||_{2}^{2}\\
&&+\left(N^{2-\mu_{1}}+1\right)||\xi||_{2}||\theta||_{2}+\left(N^{2-\mu_{1}}+N^{2-\mu_{2}}\right)||\theta||_{2}.
\end{eqnarray*}
As for the terms associated to $\widetilde{A}_{22}$ and $\widetilde{A}_{23}$, using Grisvard's lemma to bound $||{\bf u}\cdot{\bf u}||_{2}\lesssim ||{\bf u}||_{2}^{2}$, we can derive
\begin{eqnarray*}
|\sum_{i,j}(\mathcal{T}_{g}\nabla\cdot\widetilde{A}_{22},\Theta_{ij}^{N})|&\lesssim &\left(N^{2-\mu_{2}}+1\right)||\xi||_{2}||\theta||_{2},\\
|\sum_{i,j}(\mathcal{T}_{g}\nabla\cdot\widetilde{A}_{23},\Theta_{ij}^{N})|&\lesssim &\left(N^{2-\mu_{2}}+1\right)||\xi||_{2}||\theta||_{2}+||\xi||_{2}^{2}||\theta||_{2}+N^{2-\mu_{2}}||\theta||_{2}.
\end{eqnarray*}
For the control of $\mathcal{T}_{g}\nabla\cdot\widetilde{A}_{24}$, using the same arguments as before, it is not hard to derive the estimates
\begin{eqnarray*}
|\sum_{i,j}(\mathcal{T}_{g}\nabla\cdot(\xi+\sigma)P_{N}(|{\bf u}|^{2}),\Theta_{ij}^{N})|&\lesssim &\left(N^{2-\mu_{2}}+||\xi||_{2}\right)||\theta||_{2},\\
|\sum_{i,j}(\mathcal{T}_{g}\nabla\cdot(\xi+\sigma)P_{N}({\bf u}\cdot(\xi+\sigma)),\Theta_{ij}^{N})|&\lesssim &\left(N^{2-\mu_{2}}+||\xi||_{2}^{2}\right)||\theta||_{2},\\
|\sum_{i,j}(\mathcal{T}_{g}\nabla\cdot(\xi+\sigma)P_{N}(|\sigma|^{2}),\Theta_{ij}^{N})|&\lesssim &N^{2-\mu_{2}}\left(N^{2-\mu_{2}}+||\xi||_{2}\right)||\theta||_{2},\\
|\sum_{i,j}(\mathcal{T}_{g}\nabla\cdot(\xi+\sigma)P_{N}(\sigma\cdot\xi)\Theta_{ij}^{N})|&\lesssim &N^{2-\mu_{2}}\left(N^{2-\mu_{2}}+||\xi||_{2}\right)||\xi||_{2}||\theta||_{2},
\end{eqnarray*}
while for the last term we make use of (\ref{maximalt}) to have
\begin{eqnarray*}
|\sum_{i,j}(\mathcal{T}_{g}\nabla\cdot(\xi+\sigma)P_{N}(|{\xi}|^{2}),\Theta_{ij}^{N})|\lesssim \left(N^{2-\mu_{2}}+|\xi|_{\infty}\right)||\xi||_{2}^{2}||\theta||_{2}\leq \left(N^{2-\mu_{2}}+1\right)||\xi||_{2}^{2}||\theta||_{2}.
%|\sum_{i,j}(\mathcal{T}_{g}\nabla\cdot(\xi+\sigma)P_{N}({\bf u}\cdot(\xi+\sigma)),\Theta_{ij}^{N})|&\lesssim &\left(N^{2-\mu_{2}}+||\xi||_{2}^{2}\right)||\theta||_{2},\\
%|\sum_{i,j}(\mathcal{T}_{g}\nabla\cdot(\xi+\sigma)P_{N}(|\sigma|^{2}),\Theta_{ij}^{N})|&\lesssim &N^{2-\mu_{2}}\left(N^{2-\mu_{2}}+||\xi||_{2}\right)||\theta||_{2},\\
%|\sum_{i,j}(\mathcal{T}_{g}\nabla\cdot(\xi+\sigma)P_{N}(\sigma\cdot\xi)\Theta_{ij}^{N})|&\lesssim &N^{2-\mu_{2}}\left(N^{2-\mu_{2}}+||\xi||_{2}\right)||\xi||_{2}||\theta||_{2},
\end{eqnarray*}
All the previous results imply an estimate for $B_{2}$ of the form
\begin{eqnarray}
|B_{2}|&\lesssim &C(N)^{2}+||\xi||_{2}^{2}+||\theta||_{2}^{2}+||\xi||_{2}^{2}||\theta||_{2}+||\xi||_{2}||\theta||_{2}^{2},\label{BO222c}
\end{eqnarray}
where $C(N)=C(N^{3-\mu_{1}}+N^{2-\mu_{2}})$ for some constant $C$, independent of $N$. Now we define
\begin{eqnarray*}
\widetilde{Y}(t):=Y(t)+C(N)^{2}, Y(t)=\frac{\alpha\gamma}{\alpha-1}||\eta^{N}||_{2}^{2}+\frac{1}{1-\gamma}||{\bf u}^{N}||_{2}^{2},\; 0\leq t\leq t_{N}.
\end{eqnarray*}
Therefore, 
(\ref{BO222b}) and (\ref{BO222c}), applied to (\ref{BO222a}), imply
\begin{eqnarray*}
\frac{d}{dt}\widetilde{Y}(t)\leq C_{0}\widetilde{Y}(t)(1+\widetilde{Y}(t)),
\end{eqnarray*}
for some constant $C_{0}$, independent of $N$. The same argument as in Theorem \ref{th2}, cf. \cite{X}, proves the existence of some $T^{*}=T^{*}(N)>0$, given specifically as
\begin{eqnarray}
T^{*}=\frac{1}{2C_{0}}\ln\left(1+\frac{1}{\widetilde{Y}(0)}\right)=\frac{1}{2C_{0}}\ln\left(1+\frac{1}{C(N)^{2}}\right),\label{tmaximal2}
\end{eqnarray}
such that, for $0<t\leq T^{*}$
\begin{eqnarray*}
\widetilde{Y}(t)\leq  \left(1+\left(1+\frac{1}{\widetilde{Y}(0)}\right)^{1/2}\right)\widetilde{Y}(0),
\end{eqnarray*}
which, since $\widetilde{Y}(0)=C(N)^{2}$, implies
\begin{eqnarray}
{Y}(t)\leq  \sqrt{C(N)^{2}+C(N)^{4}}\lesssim C(N), \quad 0<t\leq T^{*}.\label{estY}
\end{eqnarray}
for $N$ large enough.
\end{proof}

Theorem \ref{th3} and the changes  (\ref{BO26}), (\ref{BO12c}) lead to the estimate for the original variables.
\begin{theorem}
\label{th4}
Under the conditions of Theorem \ref{th3}, there is a constant $C$, independent of $N$ such that, for $N$ large enough
\begin{eqnarray}
\max_{0\leq t\leq T^{*}}\left(||\zeta^{N}(t)-\zeta(t)||_{2}^{2}+||{\bf v}^{N}(t)-{\bf v}(t)||_{2}^{2}\right)\leq C\left(N^{3-\mu_{1}}+N^{2-\mu_{2}}\right).\label{BOest}
\end{eqnarray}
\end{theorem}
\begin{proof}
Since ${\bf v}^{N}={\bf u}^{N}, {\bf v}={\bf u}$, then we just have to study the first term on the left-hand side of (\ref{BOest}). Note that
\begin{eqnarray}
|{\bf u}^{N}|^{2}-|{\bf u}|^{2}\leq |{\bf u}^{N}-{\bf u}|(|{\bf u}^{N}|+|{\bf u}|),\label{BO223a}
\end{eqnarray}
and
\begin{eqnarray*}
||{\bf u}^{N}||_{2}\leq ||{\bf u}||_{2}+||{\bf u}^{N}-{\bf u}||_{2}.
\end{eqnarray*} 
Then, from (\ref{BO26}), (\ref{BO12c}), (\ref{BO222}), and (\ref{BO223a}), for $0\leq t\leq T^{*}$, we have
\begin{eqnarray*}
||\zeta^{N}(t)-\zeta(t)||_{2}&\leq &||\eta^{N}(t)-\eta(t)||_{2}+\frac{1}{2\gamma(1-\gamma)}\left(||P_{N}(|{\bf u}^{N}(t)|^{2}-|{\bf u}(t)|^{2})||_{2}\right.\\
&&\left.+||P_{N}(|{\bf u}(t)|^{2})-|{\bf u}(t)|^{2}||_{2}\right)\\
&\lesssim&||\eta^{N}(t)-\eta(t)||_{2}+|||{\bf u}^{N}(t)|^{2}-|{\bf u}(t)|^{2}||_{2}+N^{2-\mu_{2}}||{\bf u}(t)||_{2}^{2}\\
&\lesssim&||\eta^{N}(t)-\eta(t)||_{2}+||{\bf u}^{N}(t)-{\bf u}(t)||_{2}+N^{2-\mu_{2}}\\
&\lesssim&N^{3-\mu_{1}}+N^{2-\mu_{2}}.
\end{eqnarray*}
\end{proof}
\begin{remark}
Several comments concerning Theorem \ref{th4} are made here.
\begin{itemize}
\item From (\ref{estY}) and the hypotheses $\mu_{1}>3, \mu_{2}>2$, it holds that, for $N$ large enough, $t_{N}$ is not maximal in (\ref{maximalt}), unless $t_{N}=T^{*}$, given by (\ref{tmaximal2}).
\item The dependence of $T^{*}$ in (\ref{tmaximal2}) on $N$ implies that $T^{*}\rightarrow \infty$ as $N\rightarrow\infty$. Therefore, assuming existence of solution, globally in time, of (\ref{BO21}), Theorem \ref{th4} allows to extend the existence of numerical solution $(\zeta^{N}(t),{\bf v}^{N}(t))$ for $0\leq t\leq T^{*}$.
\item The optimality of the estimate (\ref{BOest}), in the sense given, e.~g., in \cite{DDS3} , depends on the regularity of the exact solution, represented by the parameters $\mu_{1},\mu_{2}$. Thus, we may rewrite (\ref{BOest}) in the form
\begin{eqnarray}
\left(||\zeta^{N}(t)-\zeta(t)||_{2}^{2}+||{\bf v}^{N}(t)-{\bf v}(t)||_{2}^{2}\right)\lesssim \left\{\begin{matrix}N^{3-\mu_{1}}&{\rm if}\; \mu_{1}-\mu_{2}\leq 1\\N^{2-\mu_{2}}&{\rm if}\; \mu_{1}-\mu_{2}> 1,\end{matrix}\right.\label{BOest2}
\end{eqnarray}
for $0\leq t\leq T^{*}$.
(Note that, due to the hypotheses of Theorem \ref{th4}, $\mu_{1}-\mu_{2}\geq 0$ when $\mu_{2}\leq 3$ and  $\mu_{1}-\mu_{2}< 0$ otherwise.) Therefore, (\ref{BOest2}) is optimal in $H^{2}$ for both variables in the second case, when $\mu_{1}-\mu_{2}> 1$.
\item The use of inverse inequalities (\ref{BO216}), combined with (\ref{BO215}) and  (\ref{BOest}), may lead to error estimates in higher-order spaces, under the condition of incrementing the regularity of the exact solution properly. This, if $m>2$, from  (\ref{estY}) we have, for $0<t\leq T^{*}$ and $\mu_{1}>m+1, \mu_{2}>m$,
\begin{eqnarray*}
||\theta(t)||_{m}^{2}+||\xi(t)||_{m}^{2}\lesssim N^{m-2} C(N)\lesssim N^{m+1-\mu_{1}}+N^{m-\mu_{2}},
\end{eqnarray*}
which, using (\ref{BO215}), implies
\begin{eqnarray}
\left(||\zeta^{N}(t)-\zeta(t)||_{m}^{2}+||{\bf v}^{N}(t)-{\bf v}(t)||_{m}^{2}\right)\leq C\left(N^{m+1-\mu_{1}}+N^{m-\mu_{2}}\right),
\label{BOest3}
\end{eqnarray}
for some constant $C$, independent of $N$, $0\leq t\leq T^{*}$.
\end{itemize}
\end{remark}
\section{Numerical experiments}
\label{sec3}
In this section, we complete the discretization of (\ref{BO21}) proposed here with the introduction of the time integrator and illustrate the performance of the fully discrete scheme with some numerical experiments.
\subsection{Temporal discretization and fully discrete scheme}
The semidiscrete system (\ref{BO24}) is numerically integrated in time with the Implicit Midpoint Rule (IMR) as follows: from a discretization $t_{n}=n\Delta t, n=0,\ldots,M$ of the interval $[0,T]$ and step size $\Delta t$, with $M\Delta t=T$, the fully discrete approximation at $t_{n}$ is defined as mappings $\zeta^{n}\in S_{N}, {\bf v}^{n}=(v_{1}^{n},v_{2}^{n})\in S_{N}\times S_{N}$, satisfying
\begin{eqnarray}
\frac{\zeta^{n+1}-\zeta^{n}}{\Delta t} +\frac{1}{\gamma}\mathcal{J}_{\alpha}\nabla\cdot{\bf v}^{n+1/2}-\frac{1}{\gamma}\mathcal{T}_{g}\nabla\cdot\widetilde{P}_{N}(\zeta^{n+1/2}{\bf v}^{n+1/2})&=&0,\nonumber\\
\frac{{\bf v}^{n+1}-{\bf v}^{n}}{\Delta}+(1-\gamma)\nabla \zeta^{n+1/2}-\frac{\epsilon}{2\gamma}\nabla P_{N}\left(|{\bf v}^{n+1/2}|^{2}\right)&=&0,\label{BO31}\\
\zeta^{n}\big|_{n=0}=P_{N}(\zeta_{0}),\quad {\bf v}^{n}\big|_{n=0}=\widetilde{P}_{N}({\bf v}^{0}),&&\nonumber
\end{eqnarray}
where
\begin{eqnarray*}
\zeta^{n+1/2}=\frac{\zeta^{n+1}+\zeta^{n}}{2},\quad 
{\bf v}^{n+1/2}=\frac{{\bf v}^{n+1}+{\bf v}^{n}}{2}.
\end{eqnarray*}
In terms of the Fourier coefficients $\widehat{\zeta^{n}}({\bf k}), \widehat{{\bf v}^{n}}({\bf k}), {\bf k}=(k_{x},k_{y}), -N\leq k_{x},k_{y}\leq N$, the system (\ref{BO31}) can be written as
\begin{eqnarray*}
&&\begin{pmatrix}
\widehat{\zeta^{n+1/2}}({\bf k})\\ \widehat{{\bf v}^{n+1/2}}
\end{pmatrix}=
\begin{pmatrix}\widehat{\zeta^{n}}({\bf k})\\ \widehat{{\bf v}^{n}}\end{pmatrix}+\frac{\Delta t}{2}\left(-i|{\bf k}|\mathcal{A}({\bf k})\begin{pmatrix}\widehat{\zeta^{n+1/2}}({\bf k})\\ \widehat{{\bf v}^{n+1/2}}\end{pmatrix}+\mathcal{F}^{n+1/2}\right),\nonumber\\
&&\widehat{\zeta^{n+1}}({\bf k})=2\widehat{\zeta^{n+1/2}}({\bf k})-\widehat{\zeta^{n}}({\bf k}),\;\widehat{{\bf v}^{n+1}}=2\widehat{{\bf v}^{n+1/2}}-\widehat{{\bf v}^{n}},\label{BO32}
\end{eqnarray*}
where
\begin{eqnarray}
\mathcal{A}({\bf k})&=&\begin{pmatrix}0&\frac{k_{x}}{\gamma|{\bf k}|}j_{\alpha}({\bf k})&\frac{k_{y}}{\gamma|{\bf k}|}j_{\alpha}({\bf k})\\
(1-\gamma)\frac{k_{x}}{|{\bf k}|}&0&0\\
(1-\gamma)\frac{k_{y}}{|{\bf k}|}&0&0\end{pmatrix},\label{BO33a}\\
\mathcal{F}^{n}&=&\begin{pmatrix}\frac{1}{\gamma}\frac{k_{x}\widehat{(\zeta^{n} v_{1}^{n})}({\bf k},t)+k_{y}\widehat{(\zeta^{n} v_{2}^{n})}({\bf k},t)}{1+g({\bf k})}\\
\frac{1}{2\gamma}{k_{x}}\widehat{|{\bf v}^{n}|^{2}}({\bf k},t)\\
\frac{1}{2\gamma}{k_{y}}\widehat{|{\bf v}^{n}|^{2}}({\bf k},t)\end{pmatrix}.\label{BO33b}
\end{eqnarray}
Some details on the implementation are given here. The fully discrete system (\ref{BO31}) is typically implemented in collocation form on a general interval $\Omega=(a,b)^{2}$, from a uniform grid of collocation points $z_{ij}=(x_{i},y_{j})$ with
\begin{eqnarray}
x_{i}=a+ih, y_{j}=a+jh,\; i,j=0,\ldots,2N-1,\; h=\frac{b-a}{2N},\label{BO34}
\end{eqnarray}
and the numerical approximation $\zeta^{n}, {\bf v}^{n}$ at $t_{n}$ is represented in a nodal basis of $S_{N}$ from the nodes (\ref{BO34})
\begin{eqnarray*}
\zeta^{n}(x,y)=\sum_{i,j}\zeta^{n}(z_{ij})\psi_{ij}(x,y),\quad 
{\bf v}^{n}(x,y)=\sum_{i,j}{\bf v}^{n}(z_{ij})\psi_{ij}(x,y),
\end{eqnarray*}
with 
\begin{eqnarray*}
\psi_{ij}(x,y)=\varphi_{i}(x)\phi_{j}(y),\; \varphi_{i}(x_{k})=\delta_{ik}, \phi_{j}(y_{k})=\delta_{jk}.
\end{eqnarray*}
Thus, $\zeta^{n}$ and ${\bf v}^{n}$ are represented in the collocation formulation by the nodal values at the $z_{ij}$ disposed as matrices. The temporal advance $n\mapsto n+1$ is made in Fourier space: from the values $(\zeta^{n}(z_{ij}), {\bf v}^{n}(z_{ij})), i,j=0,\ldots2N-1$, FFT techniques are used to compute the discrete Fourier coefficients $\widehat{\zeta^{n}}({\bf k}), \widehat{{\bf v}^{n}}({\bf k}), {\bf k}=(k_{x},k_{y}), -N\leq k_{x},k_{y}\leq N-1$. Then the coefficients $\widehat{\zeta^{n+1/2}}({\bf k}), \widehat{{\bf v}^{n+1/2}}({\bf k})$ are obtained iteratively from the fixed point procedure
\begin{eqnarray}
\left(I_{3}+i|{\bf k}|\frac{\Delta t}{2}\mathcal{A}({\bf k})\right)\begin{pmatrix}\zeta^{[\nu+1]}\\{\bf v}^{[\nu+1]}\end{pmatrix}=\begin{pmatrix}\widehat{\zeta^{n}}({\bf k})\\ \widehat{{\bf v}^{n}}\end{pmatrix}+\frac{\Delta t}{2}\mathcal{F}^{[\nu]},\label{BO36}
\end{eqnarray}
where the change of interval is reflected replacing ${\bf k}$ by $\frac{2\pi}{b-a}{\bf k}$. After the computation of $\widehat{\zeta^{n+1}}({\bf k}), \widehat{{\bf v}^{n+1}}({\bf k})$, we recover the values $(\zeta^{n+1}(z_{ij}), {\bf v}^{n+1}(z_{ij})), i,j=0,\ldots2N-1$, via IFFT. 

As mentioned in the introduction, the choice of IMR as time integrator is mainly due to its well-known stability properties and the typically stiff character of the semidiscrete system (\ref{BO21}). The influence of the drawback imposed by the iteration (\ref{BO36}) of the internal stage in the performance of the fully discrete scheme will be discussed in the next section.
\subsection{Numerical experiments}
We complete the presentation of the fully discrete scheme (\ref{BO31}) with some numerical experiments for the approximation of both ILW and BO systems.
\subsubsection{Experiments of validation}
We start analyzing some points of the performance of the full discretization. For $\alpha=2,\gamma=0.5$, and including the corresponding forcing term, we approximate the periodic problem (\ref{BO21}) for the ILW system on $(0,2\pi)^{2}$ for which
\begin{eqnarray}
\zeta(x,y,t)&=&e^{t}\sin(x)\cos(y),\nonumber\\
u_{1}(x,y,t)&=&e^{t}\cos^{2}(x)\sin(y),\nonumber\\
u_{2}(x,y,t)&=&e^{t}\cos(x)\sin(y),\quad x,y\in (0,2\pi),\label{BO37}
\end{eqnarray}
is solution. Taking $N=8$, in Tables \ref{BOtav1} and \ref{BOtav2}, we compute the $L^{2}$ and $L^{\infty}$ errors (denoted respectively by $E_{2}$ and $E_{\infty}$) at final time $T=1$ for the $\zeta$ and $v_{1}$ components. The regularity of (\ref{BO37}) and the error estimates (\ref{BOest}) ensure spectral convergence for the spatial discretization and the Tables show convergence in time of order two given by the IMR. 
\begin{table}[htbp]
\begin{tabular}{c|c|c|c|c|}
$\Delta t$& $E_{2}(\zeta)$&Rate&$E_{2}(v_{1})$&Rate\\
\hline
$5\times 10^{-2}$&$1.4470\times 10^{-3}$&&$9.1621\times 10^{-4}$&\\
$2.5\times 10^{-2}$&$3.9118\times 10^{-4}$&$1.887$&$2.4968\times 10^{-4}$&$1.876$\\
$1.25\times 10^{-2}$&$9.8845\times 10^{-5}$&$1.985$&$6.5391\times 10^{-5}$&$1.933$\\
$6.25\times 10^{-3}$&$2.4736\times 10^{-5}$&$1.999$&$1.6690\times 10^{-5}$&$1.970$\\
$3.125\times 10^{-3}$&$6.1846\times 10^{-6}$&$2.000$&$4.1974\times 10^{-6}$&$1.991$\\
$1.5625\times 10^{-3}$&$1.5461\times 10^{-6}$&$2.000$&$1.0513\times 10^{-6}$&$1.997$\\
\hline
\end{tabular}
\caption{$L^{2}$ errors at $T=1$ and temporal convergence rates for the method (\ref{BO31}) with respect to (\ref{BO37}) with $N=8$, $\alpha=2,\gamma=0.5$.\label{BOtav1}}
\end{table}
\begin{table}[htbp]
\begin{tabular}{c|c|c|c|c|}
$\Delta t$& $E_{\infty}(\zeta)$&Rate&$E_{\infty}(v_{1})$&Rate\\
\hline
$5\times 10^{-2}$&$3.4686\times 10^{-3}$&&$2.4663\times 10^{-3}$&\\
$2.5\times 10^{-2}$&$8.9726\times 10^{-4}$&$1.951$&$7.3987\times 10^{-4}$&$1.737$\\
$1.25\times 10^{-2}$&$2.3309\times 10^{-4}$&$1.945$&$1.6361\times 10^{-4}$&$2.177$\\
$6.25\times 10^{-3}$&$5.8766\times 10^{-5}$&$1.988$&$3.8669\times 10^{-5}$&$2.081$\\
$3.125\times 10^{-3}$&$1.4721\times 10^{-5}$&$1.997$&$9.5114\times 10^{-6}$&$2.024$\\
$1.5625\times 10^{-3}$&$3.6806\times 10^{-6}$&$2.000$&$2.3692\times 10^{-6}$&$2.005$\\
\hline
\end{tabular}
\caption{$L^{\infty}$ errors at $T=1$ and temporal convergence rates for the method (\ref{BO31}) with respect to (\ref{BO37}) with $N=8$, $\alpha=2,\gamma=0.5$.\label{BOtav2}}
\end{table}
As in the spectral approximation of other 2D models for internal waves, studied in \cite{D2026_1,D2026_2}, the main limitation seems to be related to the convergence of the iteration (\ref{BO36}). It is straightforward to see that the eigenvalues of the matrix $\mathcal{A}({\bf k})$ in (\ref{BO33a}) are $0,\pm\lambda({\bf k})$, where
\begin{eqnarray*}
\lambda({\bf k})=\left(\left(\frac{1-\gamma}{\gamma}\right) j_{\alpha}({\bf k})\right)^{1/2}.
\end{eqnarray*}
Then, the spectrum of the matrix in (\ref{BO36}) is $\{1,1\pm i|{\bf k}|\lambda({\bf k})\frac{\Delta t}{2}\}$. Since $\lambda({\bf k})$ has order zero and the nonlinear terms (\ref{BO33b}) in $\mathcal{F}$ are locally Lipschitz, we can ensure local convergence of (\ref{BO36}) for small enough $N\Delta t$.

The performance of the iteration can be illustrated in Table \ref{BO_tav2b}. For several values of $N$ and $\Delta t$, this shows the mean number $NIT(N)$ of iterations required to complete the integration of the example up to $T=1$, using a tolerance for the residual of the iteration of $10^{-9}$. For $N$ large enough, note that the convergence depends on the size of $N\Delta t$, with approximately $N\Delta t<2$ as a threshold to obtain convergence or to avoid a too slow performance.
\begin{table}[htbp]
\centering
\begin{tabular}{c|c|c|c|}
$N\Delta t$& $NIT(4)$&$NIT(8)$&$NIT(16)$\\
\hline
$4\times 10^{-1}$&$10.4$&$31.1$&$-$\\
$2\times 10^{-1}$&$7.1$&$13.7$&$14.7$\\
$1\times 10^{-1}$&$5.5$&$8.7$&$7.4$\\
$5\times 10^{-2}$&$4.2$&$6.2$&$5.2$\\
$2.5\times 10^{-2}$&$4.0$&$4.8$&$4.1$\\
$1.25\times 10^{-2}$&$4.0$&$3.8$&$3.3$\\
$6.25\times 10^{-3}$&$3.1$&$3.0$&$3.0$\\
$3.125\times 10^{-3}$&$3.0$&$2.9$&$2.5$\\
\hline
\end{tabular}
\caption{Mean number of iterations $NIT(N)$, required to complete the integration of (\ref{BO21}) with solution (\ref{BO37}) with (\ref{BO31})
up to $T=1$, for a tolerance in the residual of $10^{-9}$ in the iteration (\ref{BO36}).\label{BO_tav2b}}
\end{table}

A second example of validation is concerned with the simulation of line solitary wave solutions of the form
\begin{eqnarray}
\zeta(x,y,t)&=&\zeta(x-c_{s}t-x_{0}),\nonumber\\
v_{1}(x,y,t)&=&v(x-c_{s}t-x_{0}),\nonumber\\
v_{2}(x,y,t)&=&0,\label{BO38}
\end{eqnarray}
with $c_{s}\neq 0, x_{0}\in\mathbb{R}$, representing traveling waves in the $x$-direction with speed $c_{s}$ and permanent profiles $(\zeta(X-x_{0}),{v}(X-x_{0}))$ for the interfacial wave $\zeta$ and the first component $v_{1}$ of the velocity. The substitution of (\ref{BO38}) into (\ref{BO11b}) or (\ref{BO12b}) implies that $(\zeta(x-c_{s}t-x_{0}),v(x-c_{s}t-x_{0}))$ must be a solitary wave solution of the 1D version of the corresponding system. These solutions are not analytically known, although they can be computed numerically with high accuracy, \cite{BonaDM2021}. The profiles $(\zeta(x-x_{0}),{v}(x-x_{0}))$ are shown to decay exponentially in the ILW case and algebraically in the BO case, as in the corresponding unidirectional models, \cite{AnguloS2019}.

The simulation of (\ref{BO38}) is performed as follows. First, approximate profiles $(\zeta(x-x_{0}),{v}(x-x_{0}))$ at $x_{0}=0$ are computed with one of the procedures introduced in \cite{BonaDM2021} on a long enough interval $(-L,L), L>0$. Then the ivp is approximated by the periodic ivp (\ref{BO21}) on $(-L,L)^{2}$ with $\zeta_{0}(x,y)=\zeta(x), {\bf v}^{0}(x,y)=(v(x),0)$ as initial condition, \cite{BChen}. This is now integrated by the fully discrete scheme (\ref{BO31}) up to some final time $T$. The corresponding numerical solution at $T$ is then compared with the computed profile $(\zeta(x-x_{0}),{v}(x-x_{0}))$, centered at $x_{0}=c_{s}T$.

\begin{figure}[htbp]
\centering
%\subfigure
%{\includegraphics[width=5cm]{ISAAC2025_fig1a.eps}}
\subfigure[]
{\includegraphics[width=0.45\columnwidth]{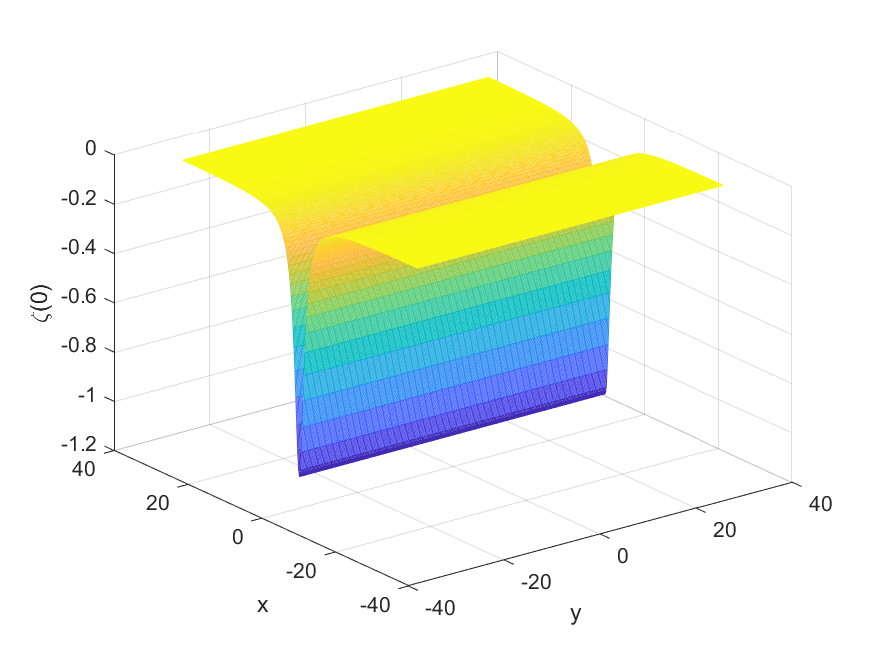}}
\subfigure[]
{\includegraphics[width=0.45\columnwidth]{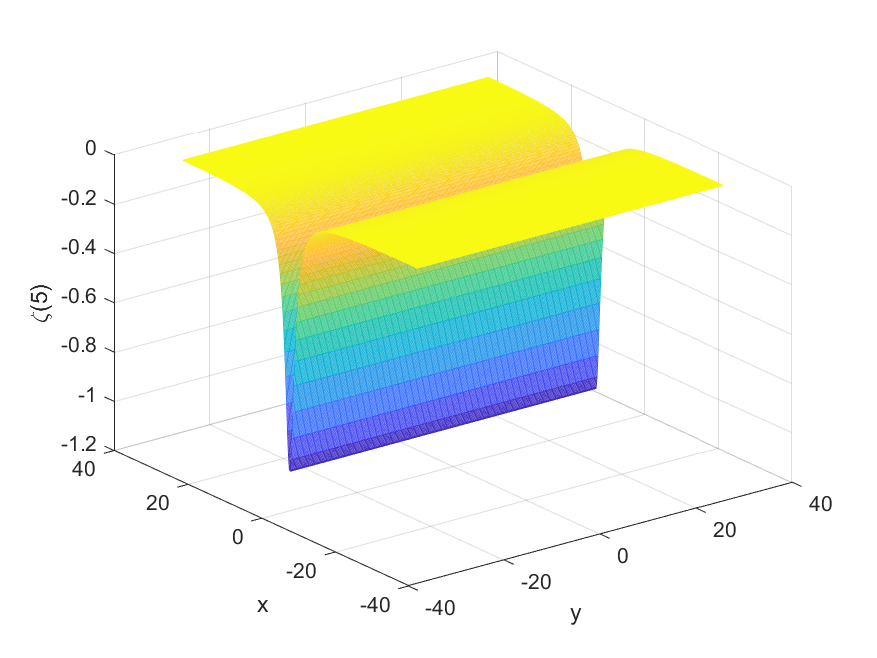}}
\subfigure[]
{\includegraphics[width=0.45\columnwidth]{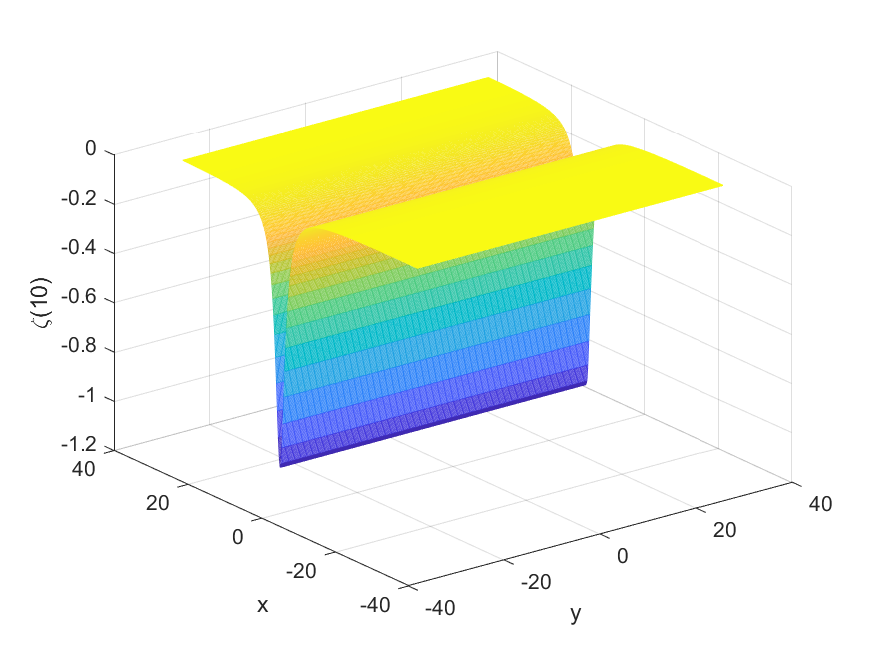}}
\subfigure[]
{\includegraphics[width=0.45\columnwidth]{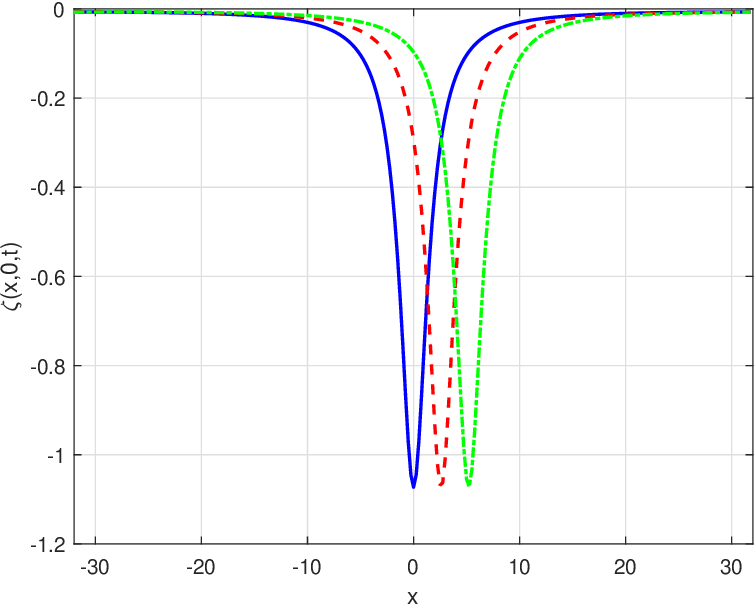}}
\caption{(a)-(c) $\zeta$ component of the numerical solution at times $t=1,5,10$; (d) Cross-sections at $y=0$.}
\label{BOFIG1}
\end{figure}

Focused on the BO system with $\alpha=1.2, \gamma=0.8, c_{s}=0.52, T=10$, and using the values $L=32, N=256$, this comparison, with the corresponding rates, is displayed, for the $L^{2}$ and $L^{\infty}$ norms and different time step sizes, in Tables \ref{BOtav3} and \ref{BOtav4}.
\begin{table}[htbp]
\begin{tabular}{c|c|c|c|c|}
$\Delta t$& $E_{2}(\zeta)$&Rate&$E_{2}(u_{2})$&Rate\\
\hline
$5\times 10^{-2}$&$3.3998\times 10^{-3}$&&$1.4809\times 10^{-3}$&\\
$2.5\times 10^{-2}$&$8.4768\times 10^{-4}$&$2.0039$&$3.6934\times 10^{-4}$&$2.0035$\\
$1.25\times 10^{-2}$&$2.0914\times 10^{-4}$&$2.190$&$9.1203\times 10^{-5}$&$2.0178$\\
$6.25\times 10^{-3}$&$4.9619\times 10^{-5}$&$2.0755$&$2.1709\times 10^{-5}$&$2.0708$\\
\hline
\end{tabular}
\caption{$L^{2}$ errors at $T=10$ and temporal convergence rates for the method (\ref{BO31}) with respect to (\ref{BO38}) with $L=32, N=256$, $\alpha=1.2,\gamma=0.8, c_{s}=0.52$.\label{BOtav3}}
\end{table}
\begin{table}[htbp]
\begin{tabular}{c|c|c|c|c|}
$\Delta t$& $E_{\infty}(\zeta)$&Rate&$E_{\infty}(u_{2})$&Rate\\
\hline
$5\times 10^{-2}$&$4.0803\times 10^{-4}$&&$1.7965\times 10^{-4}$&\\
$2.5\times 10^{-2}$&$1.0176\times 10^{-4}$&$2.0035$&$4.4810\times 10^{-5}$&$2.0033$\\
$1.25\times 10^{-2}$&$2.5146\times 10^{-5}$&$2.0167$&$1.1082\times 10^{-5}$&$2.0156$\\
$6.25\times 10^{-3}$&$5.9908\times 10^{-6}$&$2.0695$&$2.64923\times 10^{-6}$&$2.0646$\\
\hline
\end{tabular}
\caption{$L^{\infty}$ errors at $T=10$ and temporal convergence rates for the method (\ref{BO31}) with respect to (\ref{BO38}) with $L=32, N=256$, $\alpha=1.2,\gamma=0.8, c_{s}=0.52$.\label{BOtav4}}
\end{table}
The $\zeta$ component of the numerical solution is shown in Figure \ref{BOFIG1}, where the evolution of the wave of depression as localized traveling wave of permanent profile with constant speed is illustrated.

\subsubsection{Experiments on the dynamics of the models}
The performance and accuracy of the full discretization, observed in the two previous examples, can also be checked by studying computationally some aspects of the dynamics of the models (\ref{BO11b}) and (\ref{BO12b}) from different types of initial conditions. In order to approximate the ivp by some periodic ivp (as in the last example of the previous section) we  focus on the evolution of initial conditions with localized data, decaying to zero at infinity. We show results for the BO system (\ref{BO12b}) since the computations are generally harder in that case (due, for example, to the algebraic decay of the solitary waves; while for the ILW system the decay is exponential).

\begin{figure}[htbp]
\centering
%\subfigure
%{\includegraphics[width=5cm]{ISAAC2025_fig1a.eps}}
\subfigure[]
{\includegraphics[width=0.45\columnwidth]{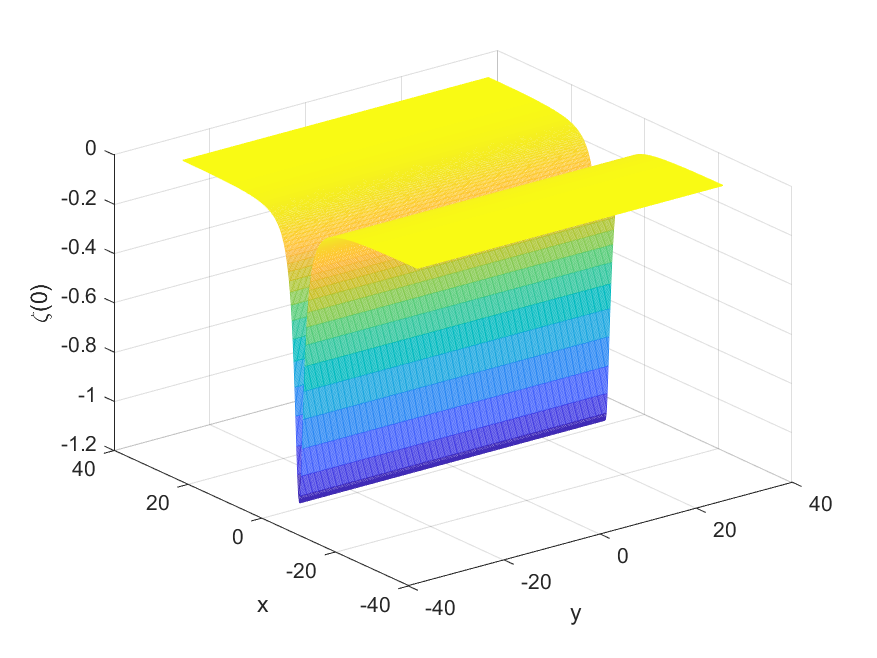}}
\subfigure[]
{\includegraphics[width=0.45\columnwidth]{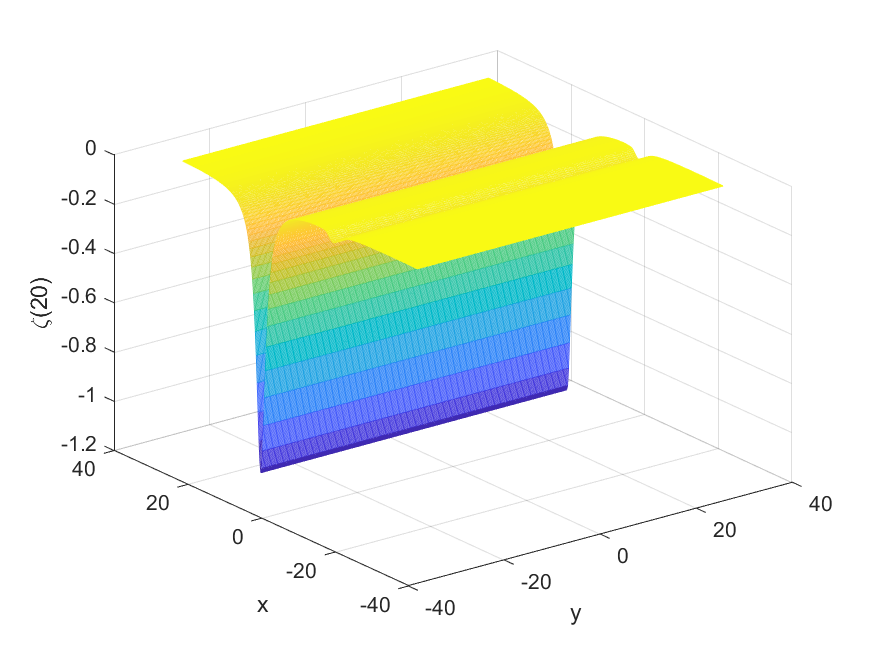}}
\subfigure[]
{\includegraphics[width=0.45\columnwidth]{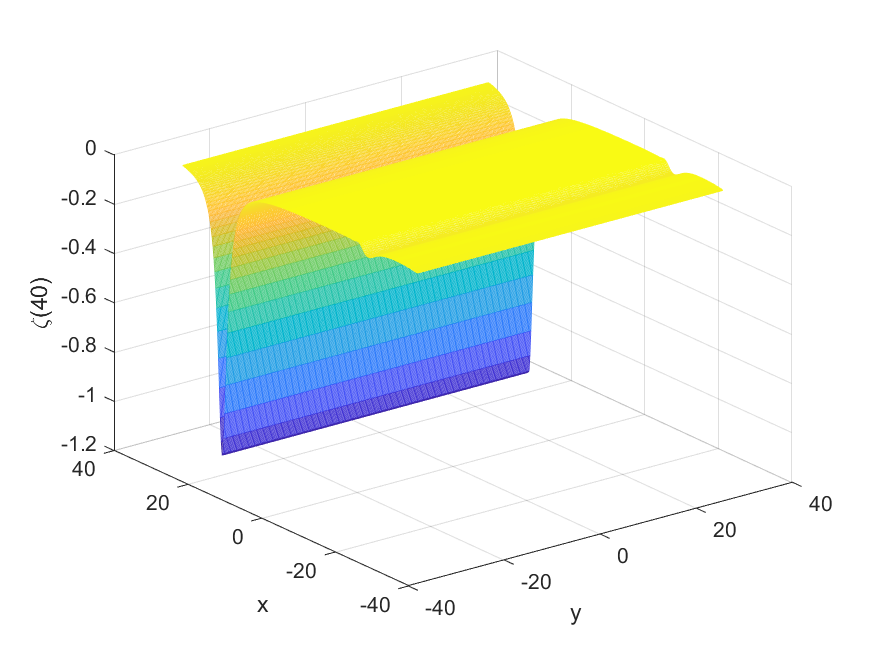}}
\subfigure[]
{\includegraphics[width=0.45\columnwidth]{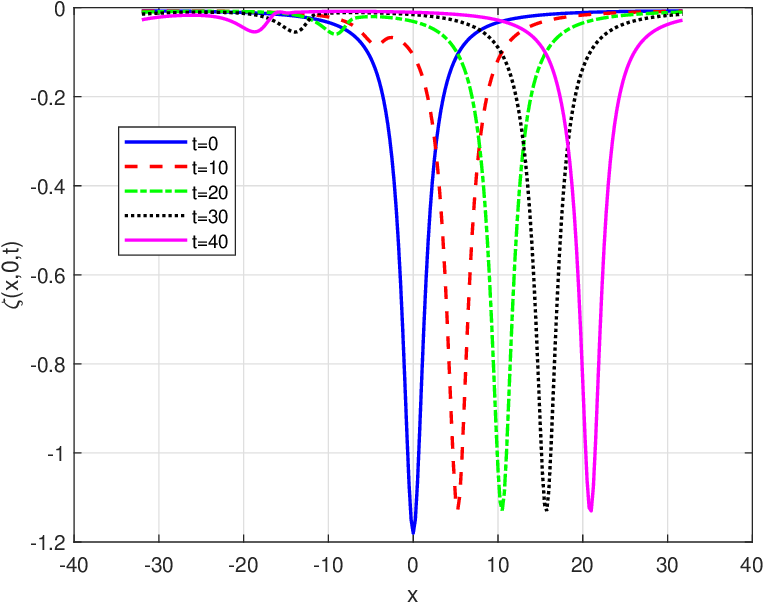}}
\caption{(a)-(c) $\zeta$ component of the numerical solution at times $t=0,20,40$ from an  initial condition (\ref{BO39}) with  $A_{1}=1.1, A_{2}=1$; (d) Cross-sections at $y=0$.}
\label{BOFIG2}
\end{figure}

By way of illustration, two types of initial waves will be considered. A first experiment is concerned with perturbations of the line solitary wave (\ref{BO38}). A typical computation here, \cite{BonaDM2022}, will consist of generating the 1D solitary wave profile $(\zeta,{v})$, consider initial conditions for our code of the perturbed form
\begin{eqnarray}
\zeta_{0}(x,y)&=&A_{1}\zeta(x),\nonumber\\
v_{1}^{0}(x,y)&=&A_{2}v(x),\nonumber\\
v_{2}(x,y,t)&=&0,\label{BO39}
\end{eqnarray}
for some constants $A_{j}, j=1,2$, and monitor the evolution of the corresponding numerical solution. Taking $A_{1}=1.1, A_{2}=1$, Figure \ref{BOFIG2} shows the evolution of the $\zeta$ component of the numerical approximation up to a final time $T=40$.

The corresponding cross-sections in Figure \ref{BOFIG2}(d) suggest that the initial, perturbed profile (\ref{BO39}) evolves to some asymptotic waveform of solitary wave type, traveling in the $x$-direction, along with a smaller, slower structure behind. According to this and other related experiments (with different vaues of $A_{j}, j=1,2$, not shown here) the general structure of the solution is kept (a modified waveform going asymptotically to some solitary wave plus a small tail behind) and the changes are just quantitative, relative to the amplitude of the elements. Furthermore,
the structure of the small tail seems to be mainly like that of a solitary wave or a nonlinear wavelet with probably some dispersive component. Out of the direction of propagation, other transverse perturbations were not observed.

\begin{figure}[htbp]
\centering
%\subfigure
%{\includegraphics[width=5cm]{ISAAC2025_fig1a.eps}}
\subfigure[]
{\includegraphics[width=0.45\columnwidth]{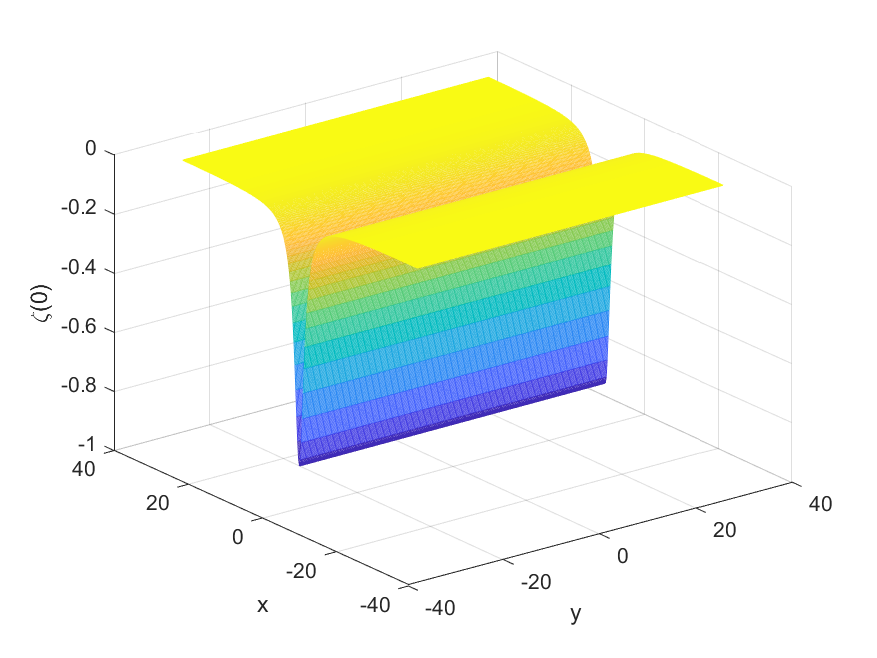}}
\subfigure[]
{\includegraphics[width=0.45\columnwidth]{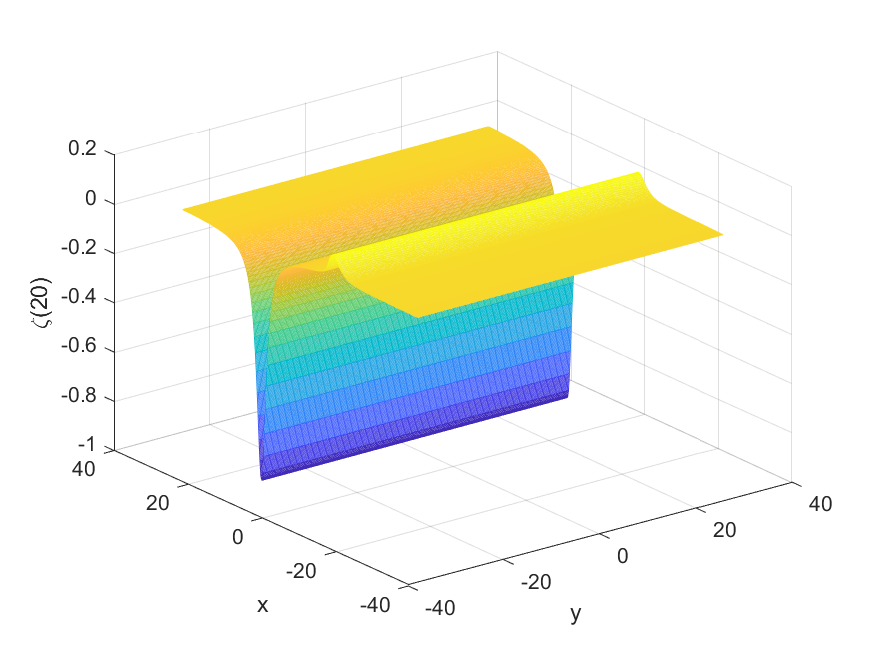}}
\subfigure[]
{\includegraphics[width=0.45\columnwidth]{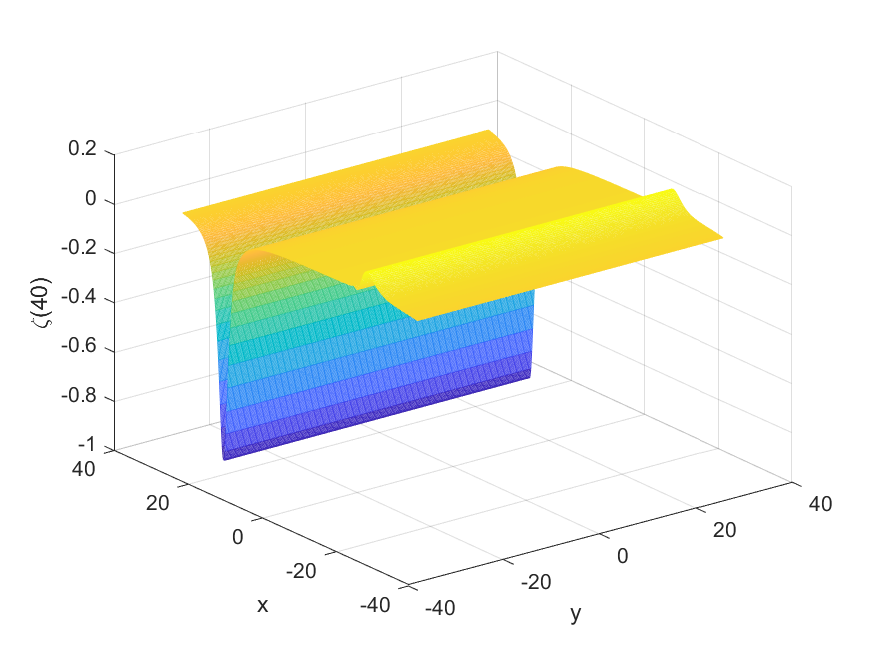}}
\subfigure[]
{\includegraphics[width=0.45\columnwidth]{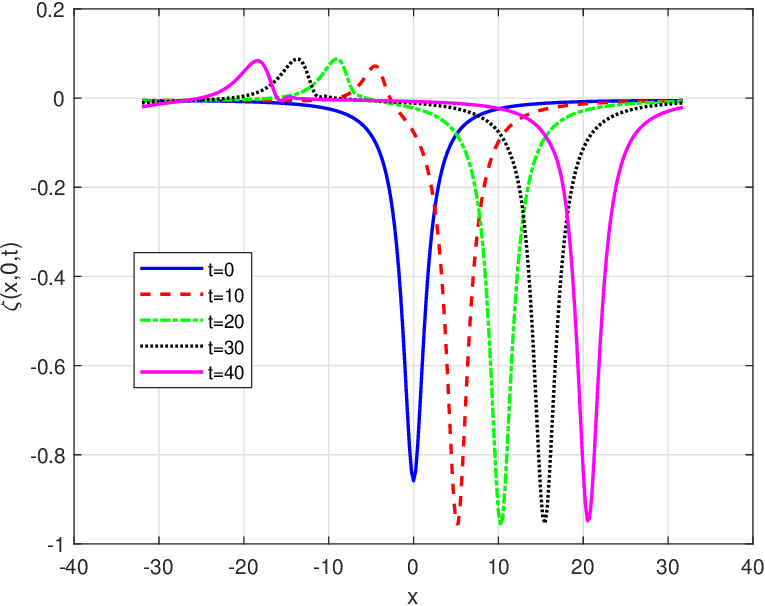}}
\caption{(a)-(c) $\zeta$ component of the numerical solution at times $t=0,20,40$ from an  initial condition (\ref{BO39}) with  $A_{1}=0.8, A_{2}=1.1$; (d) Cross-sections at $y=0$.}
\label{BOFIG2b}
\end{figure}

%\begin{figure}[htbp]
%\centering
%%\subfigure
%%{\includegraphics[width=5cm]{ISAAC2025_fig1a.eps}}
%\subfigure[]
%{\includegraphics[width=0.45\columnwidth]{BO2D_Fig1.eps}}
%\subfigure[]
%{\includegraphics[width=0.45\columnwidth]{BO2D_Fig2.eps}}
%\subfigure[]
%{\includegraphics[width=0.45\columnwidth]{BO2D_Fig3.eps}}
%\subfigure[]
%{\includegraphics[width=0.45\columnwidth]{BO2D_Fig4.eps}}
%\caption{(a)-(d) Cross-sections at $y=0$ of the $\zeta$ component of the numerical solution at times $t=10,20,30,40$ from an initial condition (\ref{BO39}) with  $A_{1}=1.1, A_{2}=1$.}
%\label{BOFIG3}
%\end{figure}
The emerging small amplitude structure may be of elevation type, as shown in Figure \ref{BOFIG2b}, corresponding to a initial condition of the form (\ref{BO39}) with  $A_{1}=0.8, A_{2}=1.1$.
%\begin{figure}[htbp]
%\centering
%%\subfigure
%%{\includegraphics[width=5cm]{ISAAC2025_fig1a.eps}}
%\subfigure[]
%{\includegraphics[width=0.45\columnwidth]{BO2D_Fig1.eps}}
%\subfigure[]
%{\includegraphics[width=0.45\columnwidth]{BO2D_Fig2.eps}}
%\subfigure[]
%{\includegraphics[width=0.45\columnwidth]{BO2D_Fig3.eps}}
%\subfigure[]
%{\includegraphics[width=0.45\columnwidth]{BO2D_Fig4.eps}}
%\caption{(a)-(d) Cross-sections at $y=0$ of the $\zeta$ component of the numerical solution at times $t=10,20,30,40$ from an initial condition (\ref{BO39}) with $A_{1}=0.8, A_{2}=1.1$.}
%\label{BOFIG3b}
%\end{figure}

A second group of experiments presented here is concerned with the presence of dispersive effects in the evolution. The study of small dispersive tails, performed in \cite{BonaDM2022} for the 1D model and with respect to perturbations of solitary waves, can be adapted and generalized to the 2D case, assuming that the linearized system of (\ref{BO13a}) 
\begin{eqnarray}
\left[Id+g(D)\right]\zeta_t+\frac{1}{\gamma}\left[Id+\frac{\alpha-1}{\alpha}g(D)\right]\nabla\cdot{\bf v}&=&0,\nonumber\\
\partial_t{\bf v}+(1-\gamma)\nabla \zeta&=&0,\label{BO311}
\end{eqnarray}
approximates the dynamics of small dispersive solutions with some accuracy. Plane wave solutions $e^{i(k_{x}x+k_{y}y-\omega({\bf k})t)}$ of (\ref{BO311}) lead to the dispersion relation
\begin{eqnarray}
\omega({\bf k})&=&\omega{\pm}({\bf k})=\pm c_{\gamma}|{\bf k}|\phi(|{\bf k}|),\; {\bf k}=(k_{x},k_{y}),\label{BO312}\\
\phi(x)&=&\left(1-\frac{\frac{\sqrt{\mu}}{\gamma}g(x)}{1+\alpha\frac{\sqrt{\mu}}{\gamma}g(x)}\right)^{1/2}, x>0;\; \phi(0)=\left\{\begin{matrix} 1&({\rm BO}\; {\rm case})\\R&({\rm ILW}\; {\rm case})\end{matrix}\right.,\nonumber\\
R&=&\left(\frac{\gamma+(\alpha-1)\frac{\sqrt{\mu}}{\mu_{2}}}{\gamma+\alpha\frac{\sqrt{\mu}}{\mu_{2}}}\right)^{1/2},\; c_{\gamma}=\sqrt{\frac{1-\gamma}{\gamma}}.\nonumber
\end{eqnarray}
Note that the dispersion relation (\ref{BO312}) is rotationally invariant. If 
\begin{eqnarray*}
v_{\pm}({\bf k})=\frac{\omega_{\pm}({\bf k})}{|{\bf k}|},\; {\bf k}\neq 0,
\end{eqnarray*}
then, according to the properties of $\phi$, proved in \cite{BonaDM2022}, it holds that
\begin{eqnarray*}
-\phi(0)c_{\gamma}<v_{-}({\bf k})<v_{+}({\bf k})<\phi(0)c_{\gamma}.
\end{eqnarray*} 
Similarly, considering the group velocity $v_{g}({\bf k})=\nabla_{{\bf k}}\omega({\bf k})$ and after some computations, we have
\begin{eqnarray*}
v_{g}({\bf k})=\pm c_{\gamma}\psi(|{\bf k}|)\frac{{\bf k}}{|{\bf k}|},\; k\neq 0,\quad 
v_{g}(0)=\pm c_{\gamma}\psi(0)\begin{pmatrix}1\\1\end{pmatrix},
\end{eqnarray*}
where
\begin{eqnarray*}
\psi(x)&=&\phi(x)+x\phi'(x),\; x>0,\\
\psi(0)&=&\phi(0).
\end{eqnarray*}
The properties of $\psi$ derived in \cite{BonaDM2022} imply that
\begin{eqnarray*}
|v_{g}({\bf k})|\leq \phi(0)c_{\gamma},
\end{eqnarray*}
and the size of the velocity group is basically controlled by the speed of sound $c_{\gamma}$, as in the 1D case.

\begin{figure}[htbp]
\centering
%\subfigure
%{\includegraphics[width=5cm]{ISAAC2025_fig1a.eps}}
\subfigure[]
{\includegraphics[width=0.45\columnwidth]{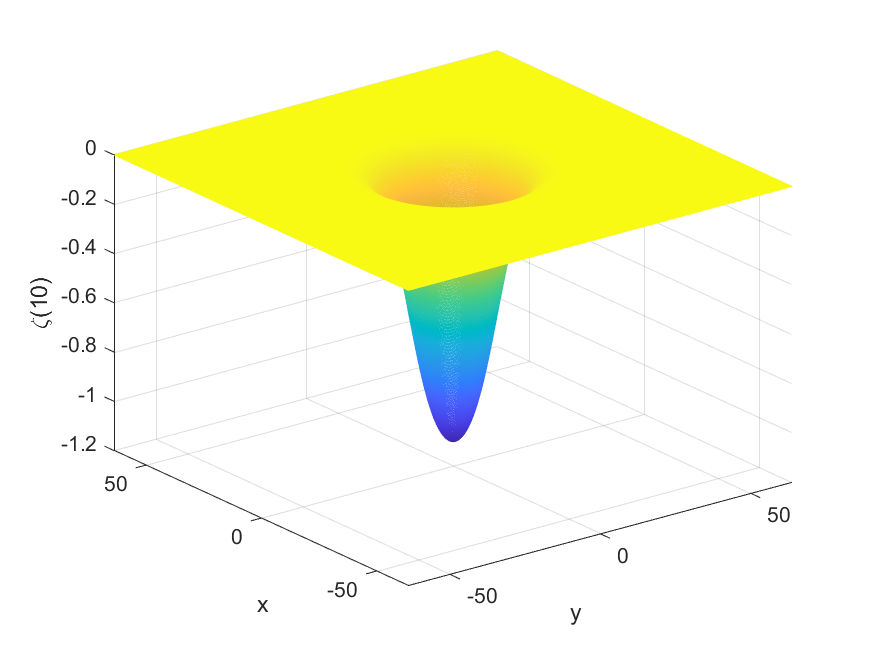}}
\subfigure[]
{\includegraphics[width=0.45\columnwidth]{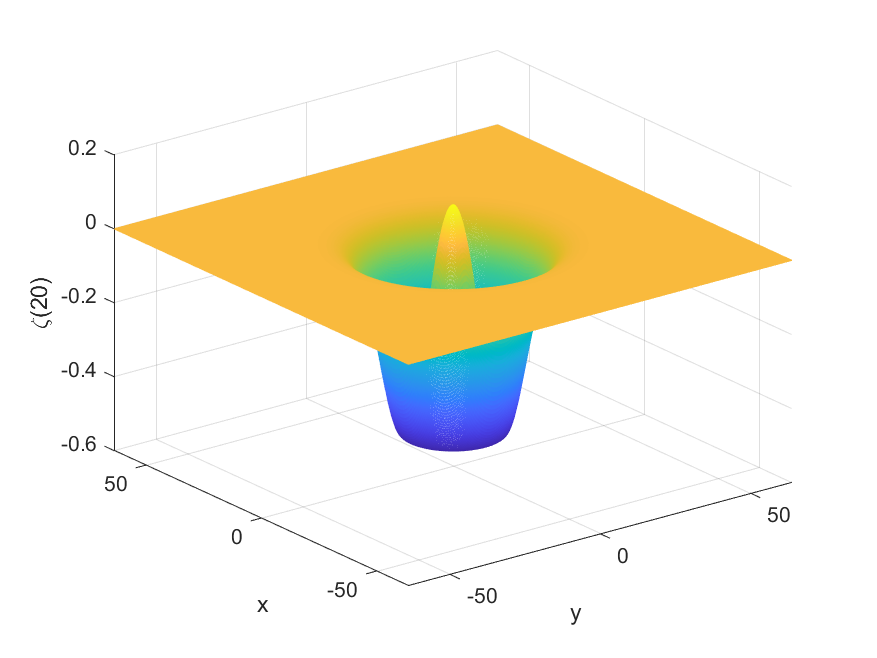}}
\subfigure[]
{\includegraphics[width=0.45\columnwidth]{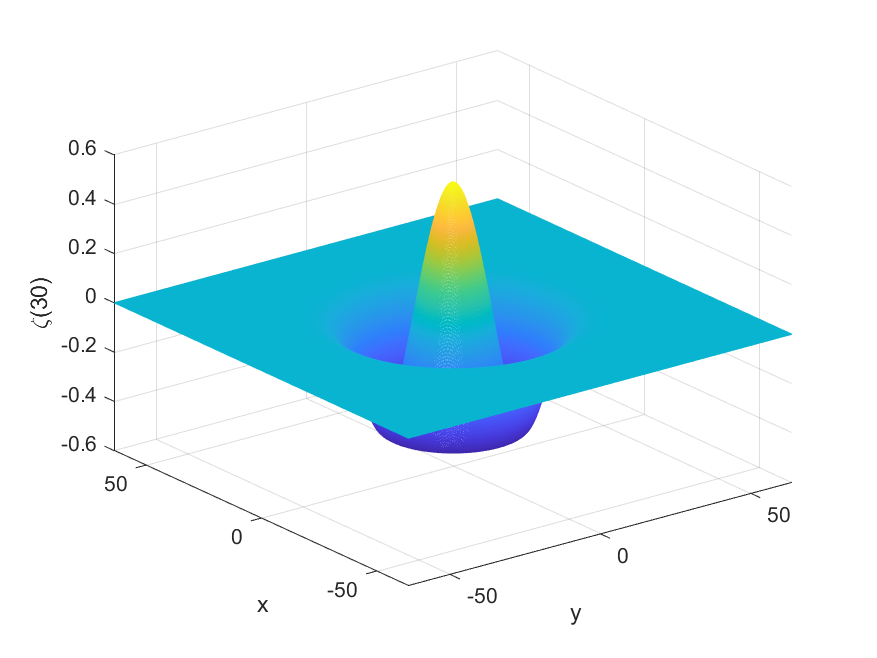}}
\subfigure[]
{\includegraphics[width=0.45\columnwidth]{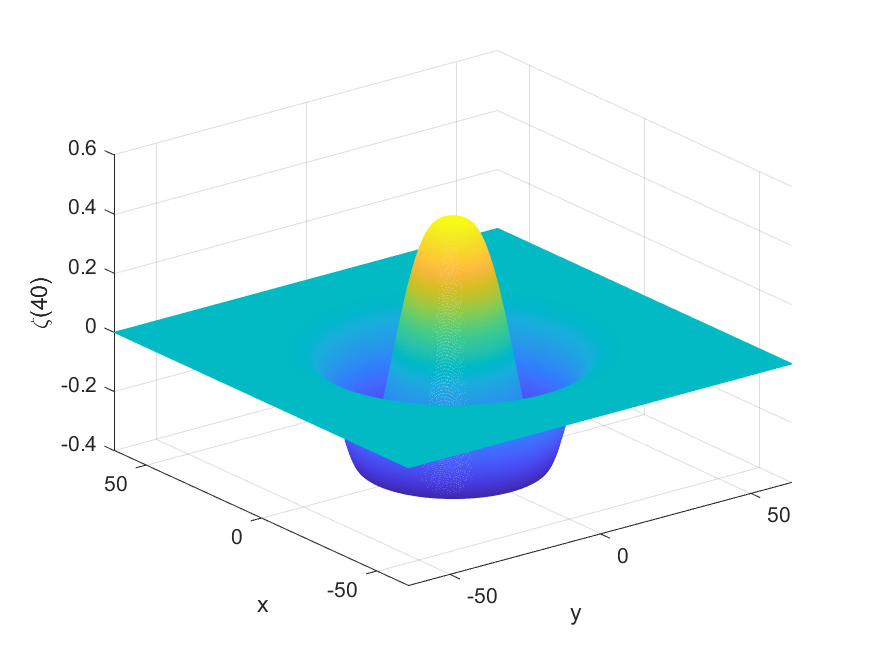}}
\caption{(a)-(d) $\zeta$ component of the numerical solution at times $t=10,20,30,40$ from an initial condition (\ref{BO310}) with  $A=2, \lambda_{1}=\lambda_{2}=10$.}
\label{BOFIG4}
\end{figure}

\begin{figure}[htbp]
\centering
%\subfigure
%{\includegraphics[width=5cm]{ISAAC2025_fig1a.eps}}
%\subfigure[]
{\includegraphics[width=0.7\textwidth]{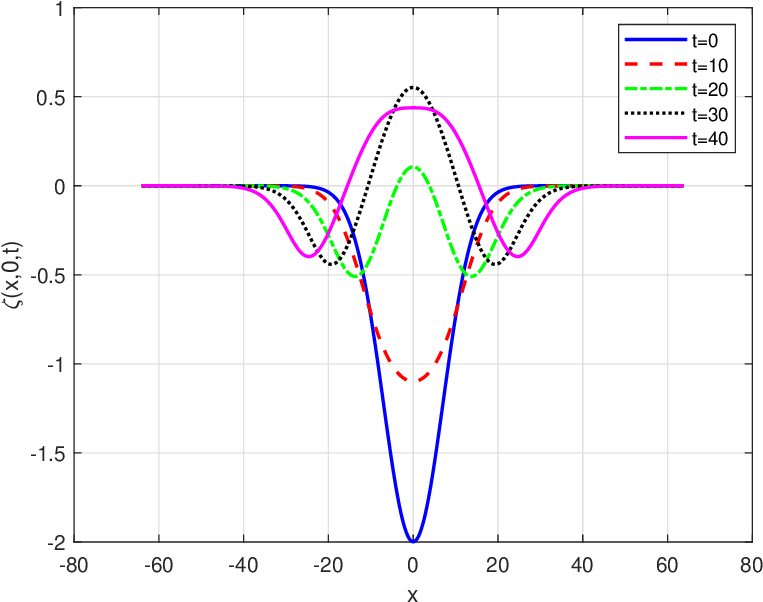}}
%\subfigure[]
%{\includegraphics[width=0.45\columnwidth]{BO2D_Fig2.eps}}
%\subfigure[]
%{\includegraphics[width=0.45\columnwidth]{BO2D_Fig3.eps}}
%\subfigure[]
%{\includegraphics[width=0.45\columnwidth]{BO2D_Fig4.eps}}
\caption{Cross-sections at $y=0$ of the $\zeta$ component of the numerical solution at times $t=0,10,20,30,40$ from an initial condition (\ref{BO310}) with  $A=2, \lambda_{1}=\lambda_{2}=10$.}
\label{BOFIG5}
\end{figure}

The previously analyzed behaviour can be observed in the evolution of the solution from initial conditions of Gaussian type
\begin{eqnarray}
\zeta_{0}(x,y)&=&Ae^{-(x/\lambda_{1})^{2}-(y/\lambda_{2})^{2}},\nonumber\\
v_{1}^{0}(x,y)&=&0),\nonumber\\
v_{2}(x,y,t)&=&0,\label{BO310}
\end{eqnarray}
for constants $A, \lambda_{j}, j=1,2$. By taking $A=2, \lambda_{1}=\lambda_{2}=10$, the dynamics of the $\zeta$ component of the numerical solution is illustrated in Figures \ref{BOFIG4}
and \ref{BOFIG5}. The initial Gaussian profile is decomposed into two groups, of apparently dispersive nature (up to the current final time of simulation), travelling to the right and to the left, and propagating in a symmetric, radial way.

\begin{figure}[htbp]
\centering
%\subfigure
%{\includegraphics[width=5cm]{ISAAC2025_fig1a.eps}}
\subfigure[]
{\includegraphics[width=0.45\columnwidth]{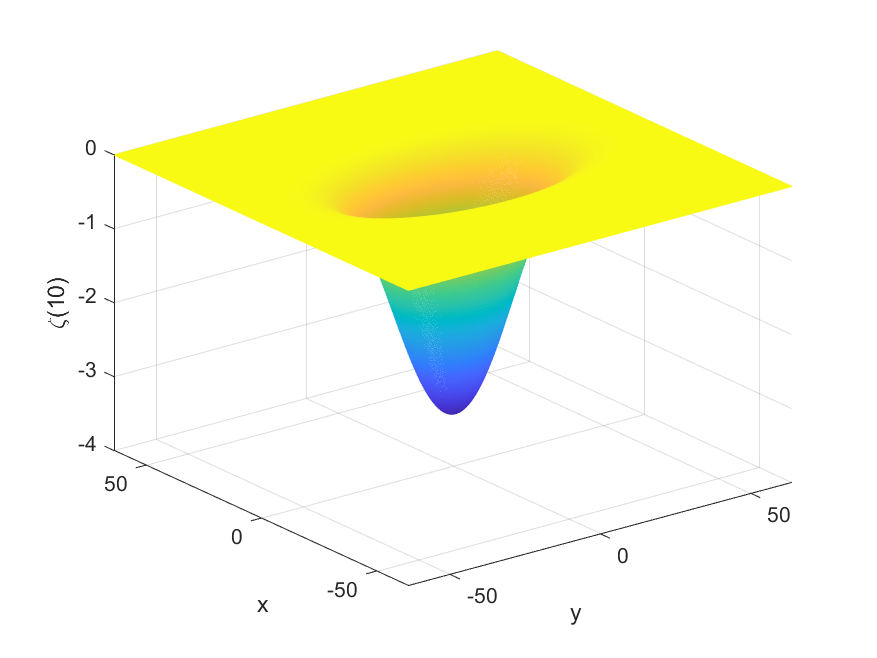}}
\subfigure[]
{\includegraphics[width=0.45\columnwidth]{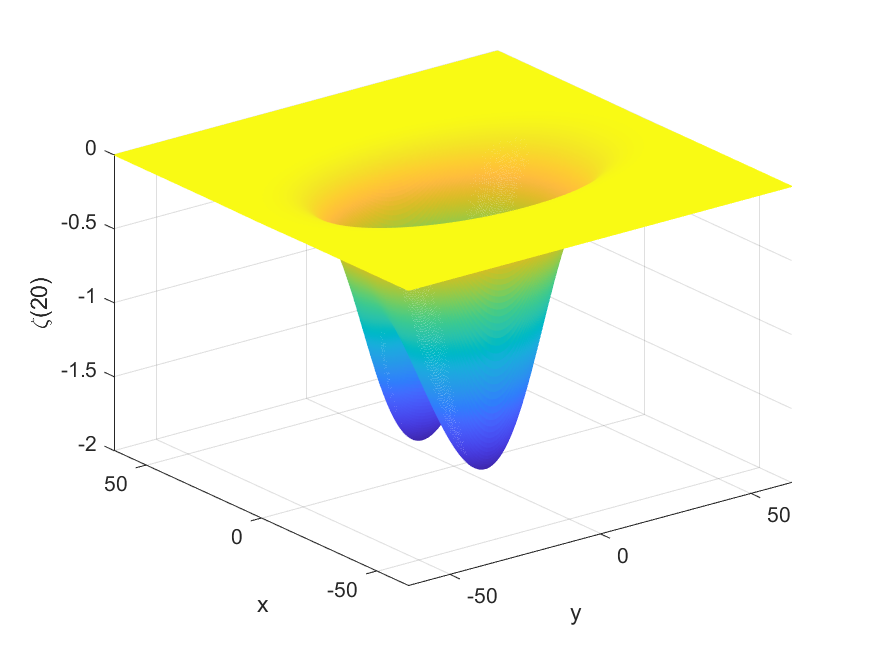}}
\subfigure[]
{\includegraphics[width=0.45\columnwidth]{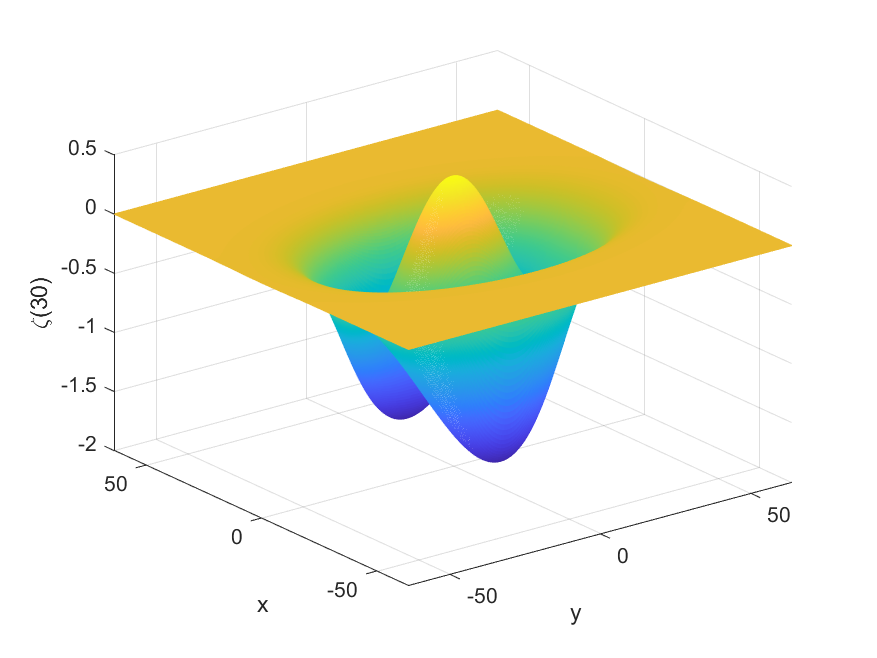}}
\subfigure[]
{\includegraphics[width=0.45\columnwidth]{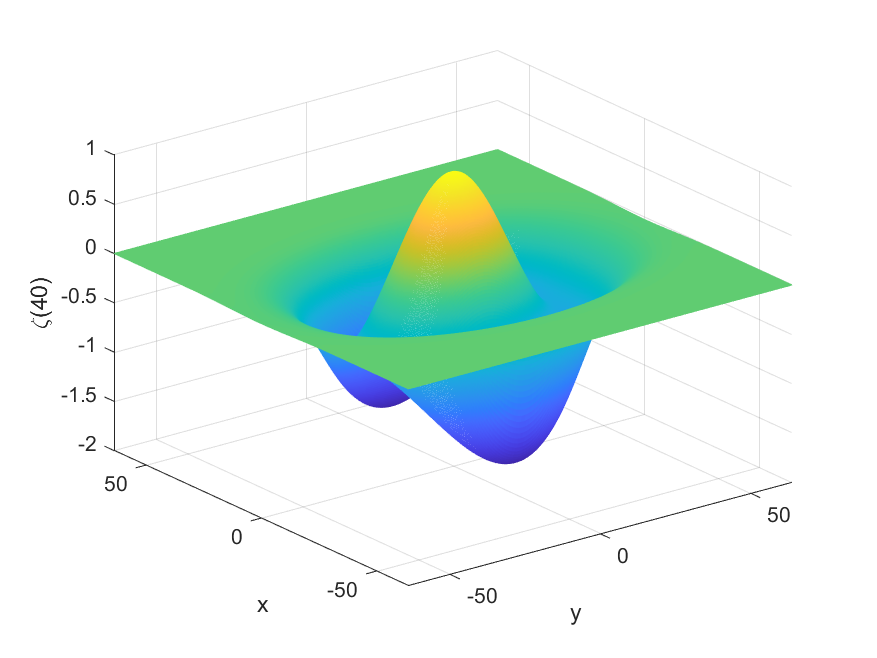}}
\caption{(a)-(d) $\zeta$ component of the numerical solution at times $t=10,20,30,40$ from an initial condition (\ref{BO310}) with  $A=5, \lambda_{1}=10, \lambda_{2}=20$.}
\label{BOFIG6}
\end{figure}

\begin{figure}[htbp]
\centering
%\subfigure
%{\includegraphics[width=5cm]{ISAAC2025_fig1a.eps}}
%\subfigure[]
{\includegraphics[width=0.7\columnwidth]{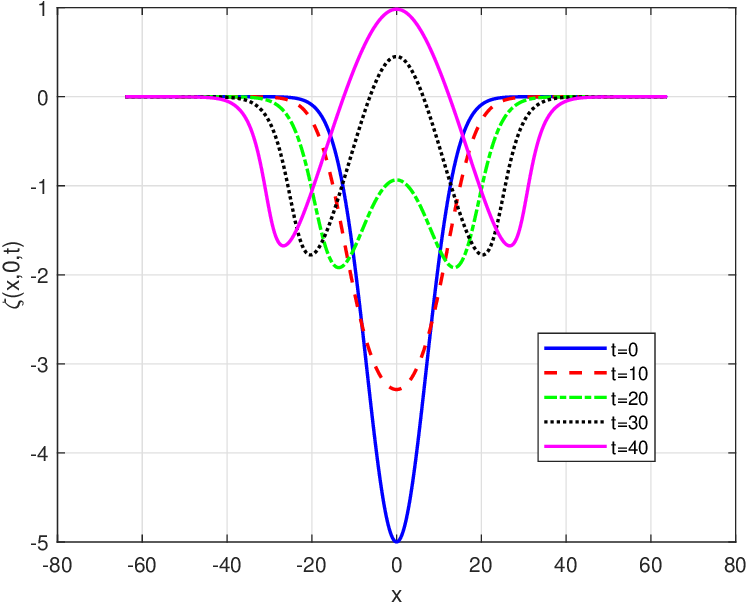}}
%\subfigure[]
%{\includegraphics[width=0.45\columnwidth]{BO2D_Fig2.eps}}
%\subfigure[]
%{\includegraphics[width=0.45\columnwidth]{BO2D_Fig3.eps}}
%\subfigure[]
%{\includegraphics[width=0.45\columnwidth]{BO2D_Fig4.eps}}
\caption{(Cross-sections at $y=0$ of the $\zeta$ component of the numerical solution at times $t=0,10,20,30,40$ from an initial condition (\ref{BO310}) with  $A=5, \lambda_{1}=10, \lambda_{2}=20$.}
\label{BOFIG7}
\end{figure}

This form of propagation does not seem to change if the initial profile is not symmetric, as Figures \ref{BOFIG6} and \ref{BOFIG7} suggest, obtained from an initial condition (\ref{BO310}) with  $A=5, \lambda_{1}=10, \lambda_{2}=20$. By increasing the amplitude of the initial wave, the asymptotic formation of solitary wave structures, traveling to the left and to the right, in a sort of resolution property and as in 1D, \cite{BonaDM2022}, is no discarded.

\section*{Acknowledgments}
This research has been supported by Ministerio de Ciencia e Innovaci\'on project PID2023-147073NB-I00, {by Junta de Castilla y Le\'on under Project VA115P25 and by Department of Education of the Junta de Castilla y Le\'on and FEDER Funds under project CLU-2025-1-02- IMUVA}.

\end{document}